\documentclass[12pt,leqno]{article}
\usepackage{amsmath,amssymb}
\newcommand{\Hom}{\mathop{\mathrm{Hom}}\nolimits}
\renewcommand{\deg}{\mathop{\mathrm{deg}}\nolimits}

\title{Some Examples of Gorenstein Liaison \\
in Codimension Three}
\author{Robin Hartshorne}

\begin{document}

\maketitle

\begin{abstract}
Gorenstein liaison seems to be the natural notion to generalize to higher
codimension the well-known results about liaison of varieties of codimension~2
in projective space.  In this paper we study points in ${\mathbb P}^3$ and
curves in ${\mathbb P}^4$ in an attempt to see how far typical codimension~2
results will extend.  While the results are satisfactory for small degree, we
find in each case examples where we cannot decide the outcome.  These examples
are candidates for counterexamples to the hoped-for extensions of codimension~2
theorems.
\end{abstract}

\vfill

\noindent
\underline{\hskip 1 in}


\clearpage

For curves in projective three-space ${\mathbb P}_k^3$, the theory of liaison,
or linkage, is classical, and is now a well-understood framework for the
classification of algebraic space curves \cite{MDP}, \cite{RHClassIII}.  This
theory has been successfully extended to schemes of codimension $2$ in any
projective space ${\mathbb P}^n$ \cite{BBM}, \cite{R}, \cite{Nollet},
\cite[Ch.~6]{M}.  Recently a number of efforts have been made to find a
suitable extension of these results in codimension
$\ge 3$ \cite{M}.  Traditional liaison uses complete intersections to link one
scheme to another.  In codimension $2$, the property of being a complete
intersection is equivalent to being arithmetically Gorenstein \cite{Serre}. 
Thus there are two natural ways to generalize.  It appears that complete
intersection liaison is too fine a relation to give analogous results in higher
codimension.  Thus attention has been focussed on Gorenstein liaison, and a
number of recent results have created an optimistic attitude that much of the
codimension $2$ case will carry over naturally to higher codimension \cite{M},
\cite{KMMNP}, \cite{CM}, \cite{N}, $\ldots$.  The purpose of this paper is to
give some examples of Gorenstein liaison for points in ${\mathbb P}^3$ and for
curves in
${\mathbb P}^4$, which suggest that the situation in codimension $\ge 3$ may be
more complicated than was initially suspected.

For points in ${\mathbb P}^3$, we show first that any set of $n$ points in
general position in a plane or on a nonsingular quadric surface can be obtained
from a single point by a sequence of ascending Gorenstein biliaisons (see \S1
below for definitions of these terms).  Thus any set of $n$ points in a plane
or on a quadric surface is glicci (in the Gorenstein liaison class of a
complete intersection).  On a nonsingular cubic surface, we can still show that
any set of $n$ points in general position is glicci, but we have to use
ascending and descending biliaisons and simple liaisons to prove this.  For a
set of $n$ points in general position in ${\mathbb P}^3$, we show that for $n
\le 19$ it is glicci, but we are unable to prove this for $n=20$.  Thus a set of
20 points in general position in ${\mathbb P}^3$ becomes a candidate for a
possible counterexample to the hope that all ACM schemes are glicci.

In ${\mathbb P}^4$, various classes of ACM curves have been shown to be glicci,
in particular, determinantal curves and ACM curves lying on general smooth ACM
rational surfaces in ${\mathbb P}^4$ \cite{KMMNP}.  We show that any general
ACM curve of degree $\le 9$ or degree $10$ and genus $6$ in ${\mathbb P}^4$ is
glicci.  Then we study ACM curves of degree $20$ and genus $26$.  There are
determinantal curves of this degree and genus, but we show that a general curve
in the irreducible component of the Hilbert scheme containing the determinantal
curves cannot be obtained by ascending Gorenstein biliaison from a line.  We do
not know if this curve is glicci, so we propose it as a candidate for an
example of an ACM curve that is not glicci.

For curves in ${\mathbb P}^4$, we consider the set of curves with Rao module
$k$, i.e., of dimension $1$ in $1$ degree only.  We call such a curve {\em
minimal} if its Rao module occurs in degree $0$.  We show that there are
minimal curves of every degree $\ge 2$.  Then we give examples that suggest
that there are curves in the liaison class of two skew lines that cannot be
reached by ascending Gorenstein biliaison from a minimal curve; and that there
are other curves with Rao module $k$ that are not in the liaison class of two
skew lines.  We will describe the examples and the evidence for these
statements below, but in most cases we cannot prove anything.

I hope that further study of these examples and others will establish whether
these guesses are correct or not, and help clarify some of the major questions
concerning Gorenstein liaison in codimension $\ge 3$.

I would like to thank Juan Migliore for his book \cite{M}, which clearly sets
out the case for Gorenstein liaison, and which stimulated this research.  I
would also like to thank him and Rosa Mir\'o--Roig and Uwe Nagel for sharing
their unpublished papers with me.  Lastly, I would like to thank the referee
for many helpful suggestions, and in particular for an idea that led to a great
improvement of Example~$4.6$.

\section{Basic Results and Questions}

Let $V_1$ and $V_2$ be two equidimensional closed subschemes without embedded
components, of the same dimension $r$ in ${\mathbb P}_k^n$, the $n$-dimensional
projective space over an algebraically closed field $k$.  We say $V_1$ and
$V_2$ are {\em linked by a complete intersection scheme $X$}, if $X$ is a
complete intersection scheme of dimension $r$ containing $V_1$ and $V_2$, and if
\begin{eqnarray*}
{\mathcal I}_{V_2,X} &\cong &\Hom({\mathcal O}_{V_1},{\mathcal O}_X), \mbox{
and} \\ 
{\mathcal I}_{V_1,X} &\cong &\Hom({\mathcal O}_{V_2},{\mathcal O}_X).
\end{eqnarray*}
Using the language of generalized divisors on Gorenstein schemes
\cite[$4.1$]{GD} we can say equivalently $V_1$ and $V_2$ are linked by a
complete intersection
$X$ if and only if there is a complete intersection scheme $S$ of dimension
$r+1$ containing $V_1$ and $V_2$, such that $X$ is an effective divisor in the
linear system $|mH|$ on $S$ for some $m > 0$, where $H$ is the hyperplane
section of $S$, and $V_2 = X - V_1$ as generalized divisors on $S$.

The equivalence relation generated by complete intersection linkage is called
{\em $CI$-liaison}.  If the equivalence can be accomplished by an even number of
links, we speak of {\em even $CI$-liaison}.

A scheme $X$ in ${\mathbb P}^n$ is called {\em arithmetically Gorenstein} (AG)
if its homogeneous coordinate ring $R/I_X$ is a Gorenstein ring, where $R =
k[x_0,\dots,x_n]$ is the homogeneous coordinate ring of ${\mathbb P}^n$, and
$I_X$ is the (saturated) homogeneous ideal of $X$.  If, in the first definition
above, we require that $X$ be an arithmetically Gorenstein scheme, instead of a
complete intersection, then we say that $V_1$ and $V_2$ are {\em linked by an
AG scheme}.  The equivalence relation generated by this kind of linkage is
called {\em Gorenstein liaison} (or {\em $G$-liaison} for short); if the
equivalence can be accomplished by an even number of $G$-links, then we speak
of {\em even Gorenstein liaison}.

One way of obtaining AG schemes is as follows.  Let $S$ be an {\em
arithmetically Cohen--Macaulay} (ACM) scheme in ${\mathbb P}^n$ (this means
that the homogeneous coordinate ring $R/I_S$ is a Cohen--Macaulay ring). 
Assume also that $S$ satisfies the property $G_1$ (Gorenstein in codimension
$1$ \cite[p.~291]{GD}), so that we can use the language of generalized
divisors.  Then any effective divisor $X$ in the linear system $|mH-K|$ on $S$,
where $m
\in {\mathbb Z}$, $H$ is the hyperplane section, and $K$ is the canonical
divisor, is an arithmetically Gorenstein scheme \cite[$4.2.8$]{M}.

Suppose now that $V_1$ and $V_2$ are divisors on an ACM scheme $S$ of dimension
$r+1$ satisfying $G_1$, let $X$ be an effective divisor in the linear system
$|mH-K|$ for some $m$, and suppose that $V_2 = X - V_1$ as generalized divisors
on $S$.  Then it is easy to see that $V_1$ and $V_2$ are linked by the AG
scheme $X$ (cf.\ proof of \cite[$4.1$]{GD} and note that since $S$
satisfies
$G_1$,
$X$ is an almost Cartier divisor on $S$ \cite[p.~301]{GD}).  In this
case we will say that $V_1$ and $V_2$ are {\em strictly $G$-linked}.  We do not
know whether the equivalence relation generated by strict $G$-linkages is
equivalent to the
$G$-liaison defined above, so we will call it {\em strict $G$-liaison}, and if
it is accomplished in an even number of steps, {\em strict even $G$-liaison}.

Combining two strict $G$-linkages gives the following result.

\bigskip
\noindent
{\bf Proposition 1.1} \cite[$5.14$]{KMMNP}, \cite[$5.2.27$]{M}. 
{\em Let
$V_1$ and
$V_2$ be effective divisors on an {\em ACM} scheme $S$ satisfying $G_1$. 
Suppose that $V_2 \in |V_1 +hH|$ for some $h \in {\mathbb Z}$.  Then $V_1$ and
$V_2$ can be strictly $G$-linked in two steps.}

\bigskip
Note that even though the statement in \cite[$5.2.27$]{M} requires $S$
smooth, the proof given there works for $S$ ACM satisfying $G_1$ if one takes
into account that any divisor $X \in |mH-K|$ is almost Cartier \cite[$2.5$]{GD}.

The proposition above motivates the following definition.  In the situation of
$(1.1)$, we say that $V_2$ is obtained by an {\em elementary Gorenstein
biliaison of height $h$} from $V_1$ \cite[$5.4.7$]{M}.  Because of the
proposition, an elementary $G$-biliaison is a strict even $G$-liaison.  If $h
\ge 0$, we call the biliaison {\em ascending}.

Now we can state some of the main questions raised by trying to generalize
codimension $2$ results to higher codimension.

\bigskip
\noindent
{\bf Question 1.2.}  a) {\em Does strict Gorenstein liaison generate the same
equivalence relation ($G$-liaison) as Gorenstein liaison?}

b) {\em Do the elementary Gorenstein biliaisons generate the same equivalence
relation as even $G$-liaison?}

\bigskip
In codimension 2, both questions reduce to $CI$-liaison, for which the answers
to parts a) and b) are both yes \cite[$4.1$, $4.4$]{GD}.

\bigskip
\noindent
{\bf Question 1.3.} a) {\em Is every {\em ACM} subscheme of ${\mathbb P}^n$ in
the
$G$-liaison class of a complete intersection (in which case we say it is {\em
glicci})?}

b) {\em Can every glicci scheme be obtained by a finite sequence of ascending
elementary $G$-biliaisons from a scheme (of the same dimension) of degree~$1$?}

\bigskip
In codimension $2$, part a) is the classical theorem of Gaeta \cite{G},
\cite{PS},
$\ldots$.  Part b) seems likely to be true, though I do not know a reference. 
Note in b) it would be equivalent to ask for ascending elementary
$G$-biliaisons starting with any complete intersection scheme.  In higher
codimension, many special cases of a) have been shown to be true
\cite{M}, \cite{KMMNP}, \cite{N}, \cite{CM}.  Closely related is the theorem of
Migliore and Nagel \cite{MN} that every ACM subscheme $X$ of ${\mathbb P}^n$
has a flat deformation to a glicci scheme, and there is also a glicci scheme
with the same Hilbert function as $X$.

\bigskip
For the following questions we limit the discussion to curves (locally
Cohen--Macaulay schemes of dimension $1$) for simplicity.  For a curve $C
\subseteq {\mathbb P}^n$ we define its {\em Rao module} to be the finite length
graded
$R$-module
$M = \oplus_{l \in {\mathbb Z}} H^1({\mathcal I}_C(l))$, where ${\mathcal I}_C$
is the ideal sheaf of $C$.  It is easy to see that even $G$-liaison preserves
the Rao module, up to shift of degrees \cite[$5.3.3$]{M}.

\bigskip
\noindent
{\bf Question 1.4.}  {\em Does the Rao module characterize the even $G$-liaison
class of a curve?  In other words, if $C$ and $C'$ are two curves with $M_{C'}
\cong M_C(h)$ for some $h \in {\mathbb Z}$, are $C$ and $C'$ in the same even
$G$-liaison class?  (In codimension $2$, this is the well-known theorem of Rao
\cite{R}.)}

\bigskip
Now we come to the problem of the structure of an even $G$-liaison class.  Let
$C \subseteq {\mathbb P}^n$ be a curve, and let ${\mathcal L}$ be the class of
all curves $C'$ in the even $G$-liaison class of $C$.  The Rao modules of
curves in ${\mathcal L}$ are all isomorphic up to shift.  As long as the Rao
module is not zero (which is equivalent to saying the curves are not ACM), one
knows that there is a minimal leftward shift of $M$ that can occur
\cite[$1.2.8$]{M}.  We denote by ${\mathcal L}_0$ the subset of ${\mathcal L}$
consisting of those curves with the leftmost possible shift of the Rao module,
and we call these {\em minimal curves}.  Let ${\mathcal L}_h$ denote the set
of curves with Rao module shifted $h$ places to the right from ${\mathcal
L}_0$, for each $h\ge 0$.  Then ${\mathcal L} = \cup {\mathcal L}_h$ for $h\ge
0$, and each one of these ${\mathcal L}_h $ for $h \ge 0$ is nonempty
\cite[$1.2.8$]{M}.

In codimension $2$ a biliaison class ${\mathcal L}$ satisfies the {\em
Lazarsfeld--Rao property} \cite[$5.4.2$]{M}.  It says that a)
${\mathcal L}_0$ is a single irreducible family of curves, and b) any curve $C
\in {\mathcal L}_h$ can be obtained by a finite sequence of ascending biliaisons
from a minimal curve, plus if necessary a deformation with constant cohomology
within the class ${\mathcal L}_h$.  (But even for curves in ${\mathbb P}^3$ it
is not known if these deformations are necessary \cite[IV, $5.4$,
p.~93]{MDP}.)  Easy examples show that in codimension $3$, ${\mathcal L}_0$ need
not consist of a single irreducible family of curves \cite[$5.4.8$]{M}.  So we
rephrase the question somewhat.

\bigskip
\noindent
{\bf Question 1.5.} a) {\em Describe the set ${\mathcal L}_0$ of minimal curves
in an even $G$-liaison class.}

b) {\em Can every curve in ${\mathcal L}_h$ for $h> 0$ be obtained from a
minimal curve by a finite sequence of ascending elementary $G$-biliaisons,
followed possibly by a flat deformation within the family ${\mathcal L}_h$?}

\bigskip
An optimist might hope for positive answers to all these questions.  However,
the examples we give below suggest that many of the answers may be no.

\section{Points in ${\mathbb P}^3$}

Any scheme of dimension zero is ACM, so in this section we will address
Question $1.3$.  The study of arbitrary zero-schemes, even in ${\mathbb P}^3$,
becomes quite complicated, so we will direct our attention to sets of reduced
points in general position.  In this section the phrase ``a general $X$ has
property $Y$'' will mean that there is a nonempty Zariski open subset of the
family of all $X$'s having the property $Y$.

We begin with points in ${\mathbb P}^2$.  In this case it is known from the
theorem of Gaeta \cite[$6.1.4$]{M} that any zero scheme in ${\mathbb P}^2$ is
{\em licci} (in the liaison class of a complete intersection), but we give a
slightly more precise statement for a general set of points and at the same
time we illustrate in a simple case the technique we will use in the later
propositions.

\bigskip
\noindent
{\bf Proposition 2.1.}  {\em A set of $n$ general points in ${\mathbb P}^2$ can
be obtained from a single point by a sequence of ascending elementary
biliaisons.}

\bigskip
\noindent
{\bf Proof.} By induction on $n$.  For $n = 1$ there is nothing to prove.  For
$n = 2$, any $2$ points lie on a line $L$.  A single biliaison of height $1$ on
$L$ reduces $2$ points to $1$ point.  For $n = 3,4,5$, a set of $n$ reduced
points, no three on a line, lies on a nonsingular conic.  These points can be
obtained by an elementary biliaison of height $1$ or $2$ from a set of $1$ or
$2$ points, and we are done by induction.

In general, let $n\ge 3$.  Then there is an integer $d \ge 2$ such that $\frac
{1}{2} (d-1)(d+2) < n \le \frac {1}{2} d(d+3)$.  Since curves of degree $d$ in
${\mathbb P}^2$ form a linear system of dimension $\frac {1}{2} d(d+3)$, any
set of $n$ points will lie on a curve of degree $d$.  Since the nonsingular
curves form a Zariski open subset of the family of all curves, a set of $n$
{\em general} points in ${\mathbb P}^2$ (in the sense mentioned above) will lie
on a nonsingular curve $C$ of degree $d$, and will form a set of $n$ general
points on $C$.  The genus of $C$ is $g = \frac {1}{2} (d-1)(d-2)$.  We use the
fact that on a nonsingular curve of genus $g$, any divisor of degree $\ge g$ is
effective.  Let $D$ be the divisor of $n$ general points on $C$, and let $H$ be
the hyperplane class on $C$.  We define a divisor $D'$ as follows.
\[
D' = \left\{ \begin{array}{ll}
D-H &\mbox{if $n = \frac {1}{2} (d-1)(d+2) + 1$} \vspace{1\jot} \\
D-2H &\mbox{if $\frac {1}{2} (d-1)(d+2) + 2 \le n \le \frac {1}{2} d(d+3)$.}
\end{array} \right.
\]
Then we verify that in either case, the degree of $D'$ is $\ge g$, so that the
divisor $D'$ is effective, and secondly that $n' = \deg D' \le \frac {1}{2}
(d-1)(d+2)$.

Now, by induction on $d$, a general set of $n'$ points can be obtained from a
single point by ascending biliaisons.  Since any $D$ as above bilinks down to a
$D'$, it follows that by bilinking up $n'$ general points on $C$, we obtain $n$
general points on $C$, as required.

\bigskip
\noindent
{\bf Proposition 2.2.}  {\em A set of $n$ general points on a (fixed)
nonsingular quadric surface $Q \subseteq {\mathbb P}^3$ can be obtained from a
single point by a finite number of ascending elementary $G$-biliaisons on $Q$.}

\bigskip
\noindent
{\bf Proof.} The method is analogous to the proof of $(2.1)$, except that now
we use both types of ACM curves on $Q$.  For $n=1$, there is nothing to prove. 
For $n = 2$, we put $2$ points on a twisted cubic curve, and then move them by
linear equivalence (a biliaison of height $0$) until they lie on a line on
$Q$.  On that line, we obtain $2$ points by a biliaison of height $1$ from $1$
point.  For $n = 3$, the points lie on a conic, and come from $1$ point by
biliaison.  For $n = 4,5$, the points lie on a twisted cubic curve, and reduce
by biliaison to $1$ or $2$ points.

Now suppose $n\ge 6$.  Then there is an integer $a \ge 2$ such that either

\begin{itemize}
\item[i)] $a^2 + a \le n \le a^2 + 2a$, or
\item[ii)] $a^2 + 2a + 1 \le n \le a^2 + 3a + 1$.
\end{itemize}
In case i) we consider the complete intersection curve $C$ of bidegree $(a,a)$
on $Q$.  It has degree $2a$ and genus $g = (a-1)^2$ and moves in a linear
system of dimension $a^2 + 2a$.  Hence $n$ general points lie on a smooth such
curve $C$, forming a divisor $D$.  The divisor $D' = D - H$ on $C$ has degree
$n' = n-2a$.  Since $n' \ge a^2 - a > g$, the divisor $D'$ is effective.  On
the other hand $n' \le a^2$, so $n'$ falls in the range i) for $a-1$ unless $n'
= a^2$, in which case it falls in range ii) for $a-1$.  By induction on $a$, a
set of $n'$ general points can be obtained by ascending elementary
$G$-biliaisons from a point, so also can $n$ points.

In case ii), we consider the ACM curve $C$ of bidegree $(a,a+1)$ on $Q$.  It
has degree $d = 2a+1$, genus $g = a(a-1)$, and moves in a linear system of
dimension $a^2 + 3a + 1$ on $Q$.  So $n$ general points form a divisor $D$ on a
nonsingular such curve $C$.  The divisor $D' = D-H$ has degree $n' = n - 2a - 1
\ge a^2 > g$, so $D'$ is effective.  On the other hand $a^2 \le n' \le a^2 + a
+ 1$, which is range ii) for $a-1$.  So by induction on $a$ again, we can
obtain $D'$ by biliaisons from a point, and $D$ by a single elementary
$G$-biliaison for $D'$ on $C$.

\bigskip
\noindent
{\bf Corollary 2.3.} {\em A set of $n$ general points on a nonsingular quadric
surface $Q$ in ${\mathbb P}^3$ is in the strict even $G$-liaison class of a
point.  In particular, it is glicci.}

\bigskip
\noindent
{\bf Proposition 2.4.} {\em A set of $n$ general points on a (fixed)
nonsingular cubic surface $S$ in ${\mathbb P}^3$ is in the same strict
Gorenstein liaison class as a point on $S$.}

\bigskip
\noindent
{\bf Proof.} A curve $C$ of degree $d$ and genus $g$ on $S$ moves in a linear
system of dimension $d+g-1$.  If a set of $n$ general points is to form a
divisor $D$ on $C$, we need $n\le d+g-1$.  In that case the linear system $D-H$
on $C$ has degree $\le g-1$, and hence may not be effective.  Thus we cannot
use Gorenstein biliaisons on the cubic surface.  Instead, we will use strict
Gorenstein liaison by AG divisors in the linear systems $|mH-K|$ on ACM curves
$C$ on $S$.

There are four types of smooth ACM curves on $S$, obtained by biliaison (of
curves) on $S$ from the line, the conic, the twisted cubic, and the hyperplane
class $H$, which is a plane cubic curve.  For $a \ge 1$ the four types are

\begin{itemize}
\item[i)] $d = 3a - 2$, $g = \frac {1}{2} (3a^2 - 7a + 4)$
\item[ii)] $d = 3a - 1$, $g = \frac {1}{2} (3a^2 - 5a + 2)$
\item[iii)] $d = 3a$, $g = \frac {1}{2} (3a^2-3a)$
\item[iv)] $d = 3a$, $g = \frac {1}{2} (3a^2 - 3a + 2)$.
\end{itemize}
In each case, one of these curves $C$ with degree $d$ and genus $g$ (we say
{\em type} $(d,g)$) moves in a linear system of dimension $d+g-1$.  On an ACM
curve $C$ of type $(d,g)$, we will consider only divisors of degree $n$, where
$g \le n \le d+g-1$.  If $n$ and $n'$ are both in this range, and if $n+n' =
\deg(mH-K)$ for some $m$, then a strict Gorenstein liaison by AG divisors in
the linear system $|mH-K|$ will transform general divisors of degree $n$ to
general divisors of degree $n'$ and vice versa.  To explain this in more
detail, let $Z$ be a set of $n$ general points on $S$.  If $n \le d + g - 1$,
then $Z$ is contained in a curve $C$ as above.  If $n' = \deg(mH-K) - n \ge
g$, then there is an effective divisor $Z'$ of degree $n'$ such that $Z + Z'
\in |mH-K|$.  Thus $Z$ and $Z'$ are linked.  The same arguments work in
reverse, starting with $Z'$, assuming $n' \le d + g - 1$ and $n \ge g$.  Hence
there are Zariski open subsets $U$ (resp.\ $U'$) of the set of all subsets of
$n$ (resp.\ $n'$) points of $S$ such that each $Z \in U$ is linked to a $Z'
\in U'$ and vice versa.  If one of these is already known to be in the strict
$G$-liaison class of a point, we conclude so is the other.  We will write $n
\leftrightarrow n'$ by
$mH - K$ on
$(d,g)$.

For $n \le 8$, we use the following liaisons.

\begin{itemize}
\item[1)] $1 \leftrightarrow 3$ by $H-K$ on $(4,1)$
\item[2)] $2 \leftrightarrow 6$ and $3 \leftrightarrow 5$ by $2H-K$ on $(5,2)$
\item[3)] $6 \leftrightarrow 8$ by $3H - K$ on $(6,3)$
\item[4)] $4 \leftrightarrow 8$ and $5 \leftrightarrow 7$ by $3H-K$ on $(6,4)$
\item[5)] $6 \leftrightarrow 7$ by $3H-K$ on $(7,5)$.
\end{itemize}
These liaisons show that any set of $n\le 8$ general points on $S$ is in the
strict $G$-liaison class of a point.  Note that we must use ascending {\em and}
descending liaisons to accomplish this.  For example, if $n = 2$, the links go
$2 \rightarrow 6 \rightarrow 7 \rightarrow 5 \rightarrow 3 \rightarrow 1$.

For $9 \le n \le 17$, we use the following links.

\begin{itemize}
\item[1)] $9 \leftrightarrow 11$ by $4H-K$ on $(7,5)$
\item[2)] $7 \leftrightarrow 13$ and $8 \leftrightarrow 12$ by $4H-K$ on $(8,7)$
\item[3)] $12 \leftrightarrow 17$ and $13 \leftrightarrow 16$ by $5H-K$ on
$(9,9)$
\item[4)] $10 \leftrightarrow 17$; $11 \leftrightarrow 16$; $12 \leftrightarrow
15$; and $13 \leftrightarrow 14$ by $5H-K$ on $(9,10)$.
\end{itemize}
Using these links a general set of $9 \le n \le 17$ points is linked down to a
set of $7$ or $8$ points treated above.

For $n \ge 18$, we find an integer $a$ such that $n_0 = \frac {3}{2} a(a-1) \le
n < n_1 = \frac {3}{2} (a+1)a$.  We divide these $n$'s into six ranges.

\begin{itemize}
\item[A)] $n_0 \le n \le n_0+2$
\item[B)] $n_0+3 \le n \le n_0 + a-1$
\item[C)] $n_0 + a$, $n_0 + a + 1$
\item[D)] $n_0 + a + 2 \le n \le n_0 + 2a - 2$
\item[E)] $n_0 + 2a - 1$, $n_0 + 2a$
\item[F)] $n_0 + 2a + 1 \le n \le n_0 + 3a - 1$.
\end{itemize}

In range $A$, we do

\begin{itemize}
\item[1)] $n_0 \leftrightarrow n_0 + 2$ by $(2a-2)H-K$ on type (i) curve.
\item[2)] $n_0 + 1 \leftrightarrow n_0 + 3a-1$ and $n_0 + 2 \leftrightarrow n_0
+ 3a - 2$ by $(2a-1)H-K$ on type (iv).
\end{itemize}

In range $B$, we do

\begin{itemize}
\item[3)] $n_0 + t \leftrightarrow n' < n_0$ by $(2a-2)H-K$ on (i).
\end{itemize}

In range $C$, we do

\begin{itemize}
\item[4)] $n_0 + t \leftrightarrow n' < n_0$ by $(2a-2)H-K$ on (ii).
\end{itemize}

In range $D$, we do links by $(2a-1)H-K$ on types (ii) and (iv).  If $a$ is
odd, $a = 2k+1$, start with $m = n_0 +3k +2$, link on type (iv).  Then
alternate linkages on type (ii) and (iv).  This covers all values of $n$ in
range $D$, starting in the middle and spiraling outward, until finally we land
in range $E$.  If $a$ is even, $a = 2k$, start with $m = n_0+3k$, and do a link
on type (ii) first, then alternate (iv) and (ii).

In range $E$, link by $(2a-1)H-K$ on type (iv), to land in range $C$.

In range $F$, link by $(2a-1)H-K$ on type (iii) to land in range $B$ or $C$.

In summary, ranges $B$ and $C$ link down to $n' < n_0$ so are ok by induction. 
Ranges $E$ and $F$ link down to ranges $B$ and $C$.  Range $D$ spirals up and
down until it lands in range $E$; and finally range $A$ links up to range $F$. 
So for example, if $n = 18$ (range $A$), the links go $18 \rightarrow 20
\rightarrow 28 \rightarrow 22 \rightarrow 16 \rightarrow 13 \rightarrow 7$,
which we did earlier.  If $n = 54$ (range $D$), the links go $54 \rightarrow 55
\rightarrow 53 \rightarrow 56 \rightarrow 52 \rightarrow 40 \rightarrow 35
\rightarrow 27 \rightarrow 23 \rightarrow 15 \rightarrow 12 \rightarrow 8
\rightarrow 6$, treated above.

\bigskip
\noindent
{\bf Corollary 2.5.} {\em A set of $n$ general points on a smooth cubic surface
in ${\mathbb P}^3$ is glicci.}

\bigskip
\noindent
{\bf Corollary 2.6.} {\em A set of $n \le 19$ general points in ${\mathbb P}^3$
is glicci.}

\bigskip
\noindent
{\bf Proof.} Indeed, since the cubic surfaces in ${\mathbb P}^3$ form a linear
system of dimension $19$, a set of $n \le 19$ general points lie on a smooth
cubic surface, and we can apply $(2.4)$.  Using $(2.2)$ we can also see that a
set of $n\le 9$ general points can be obtained from a single point by ascending
elementary $G$-biliaisons.  However, we can get a stronger result, and another
proof of $(2.6)$ by another method.

\bigskip
\noindent
{\bf Proposition 2.7.} {\em A set of $n \le 19$ general points in ${\mathbb
P}^3$ is in the strict Gorenstein liaison class of a point.  Furthermore,
if $n \ne 17,19$, it can be obtained by a sequence of ascending
elementary $G$-biliaisons from a point.}

\bigskip
\noindent
{\bf Proof.} Again we use ACM curves $(d,g)$ in ${\mathbb P}^3$, but now we
need to know how many general points of ${\mathbb P}^3$ lie on such a curve. 
Call this number $m(d,g)$.  This is no longer an elementary question, because
the families of these curves form a Hilbert scheme, not a linear system.  The
question was studied in Perrin's thesis \cite{P}, and depends on semi-stability
of the normal bundle.  Here are his results, for the ACM curves we need:

\begin{center}
\begin{tabular}{crl}
$(d,g)$ & $m$ & reference in \cite{P} \\ \hline
(1,0) & 2 \\
(2,0) & 3 \\
(3,0) & 6 \\
(4,1) & 8 \\
(5,2) & 9 \\
(6,3) & 12 & \hskip 0.15 in p.~66 \\
(7,5) & 14 & \hskip 0.15 in p.~87 \\
(8,7) & 16 &\hskip 0.15 in  p.~10; p.~116 \\
(9,9) & 18 & \hskip 0.15 in p.~87 \\
(10,11) & 20 & \hskip 0.15 in p.~66
\end{tabular}
\end{center}
In this table $m$ has the naive value $2d$, except for $(2,0)$, a conic, which
lies in a plane, so can pass through at most $3$ general points, and $(5,2)$,
which lies on a quadric surface, so can pass through at most $9$ general points.

To prove our result, we use the following biliaisons and liaisons.

\begin{itemize}
\item[1)] $1 \leftrightarrow 2$ biliaison on $(1,0)$
\item[2)] $1 \leftrightarrow 3$ biliaison on $(2,0)$
\item[3)] $1,2,3 \leftrightarrow 4,5,6$ biliaison on $(3,0)$
\item[4)] $3,4 \leftrightarrow 7,8$ biliaison on $(4,1)$
\item[5)] $4 \leftrightarrow 9$ biliaison on $(5,2)$
\item[6)] $4,5,6 \leftrightarrow 10,11,12$ biliaison on $(6,3)$
\item[7)] $6,7 \leftrightarrow 13,14$ biliaison on $(7,5)$
\item[8)] $8,9 \leftrightarrow 15,16$ biliaison on $(8,7)$
\item[9)] $9 \leftrightarrow 18$ biliaison on $(9,9)$
\item[10)] $12 \leftrightarrow 17$ liaison by $5H-K$ on $(9,9)$
\item[11)] $11 \leftrightarrow 19$ liaison by $5H-K$ on $(10,11)$.
\end{itemize}

\bigskip
\noindent
{\bf Remark 2.8.} If we consider a set of 20 general points in ${\mathbb P}^3$,
none of the above methods works.  They do not lie on a cubic surface, so we
cannot apply $(2.4)$.  They do form a divisor $D$ on an ACM curve $(10,11)$,
but $D-H$ has degree $10$, less than the genus, so it may not be effective. 
Liaison by $5H-K$ would give a divisor of degree $10$, which may not be
effective.  Liaison by $6H-K$ gives another general divisor of degree $20$, so
we get nowhere.

It is conceivable that some upward liaison may eventually lead to a zero-scheme
that can then be linked back down to a point.  Or perhaps there are other AG
schemes in ${\mathbb P}^3$ besides the ones of the form $mH-K$ on ACM curves
that we have been using.

On the other hand, it may simply be that 20 general points in ${\mathbb P}^3$
are not in the $G$-liaison class of a point, so we propose this as a potential
counterexample to Question $1.3$.  If this is so, the cone over these 20 points
would be an ACM curve in ${\mathbb P}^4$ that is not glicci.

\bigskip
\noindent
{\bf Remark 2.9.} In \cite[$3.1$]{KMMNP}, the authors prove that any standard
determinantal scheme is glicci.  Taking the $t \times t$ minors of a $t \times
(t+2)$ matrix of linear forms in ${\mathbb P}^3$ gives a zero-dimensional
determinantal scheme of degree $\frac {1}{6} (t+2)(t+1)t$, which is glicci by
the above result.  For $t = 1,2$, any set of $1$ or $4$ points in ${\mathbb
P}^3$ is determinantal.  However, for $t = 3,4$, the dimension calculation in
\cite[$10.3$]{KMMNP} show that determinantal sets of $10$ points have
codimension $3$ in zero-schemes of degree $10$, and determinantal sets of $20$
points have codimension $15$ among zero-schemes of degree $20$.

\section{ACM curves in ${\mathbb P}^4$}

In the literature, a number of special cases of ACM curves have been shown to
be glicci \cite{CM}, \cite{KMMNP}.  In this section we begin a systematic study
of ACM curves of small degree in ${\mathbb P}^4$.  We show that any general ACM
curve of degree $\le 9$, or a general ACM curve of degree $10$ and genus $6$
can be obtained by ascending Gorenstein biliaisons from a line.  On the other
hand, we show that there is an irreducible component of the Hilbert scheme of
curves of degree $20$ and genus $26$ whose general member is a smooth ACM curve
that cannot be obtained by ascending Gorenstein biliaisons from a line.  We
propose this curve as a candidate for a possible counterexample to the question
whether every ACM curve is glicci.

We start by finding a lower bound on the genus of an ACM curve in ${\mathbb
P}^4$.

\bigskip
\noindent
{\bf Proposition 3.1.}  {\em Let $C$ be a nondegenerate (i.e., not contained in
a hyperplane) {\em ACM} curve in ${\mathbb P}^4$, of degree $d$ and arithmetic
genus
$g$.  Then $d \ge 4$ and $g \ge G_{\mbox{\em min}}(d)$, where}
\[
G_{\mbox{min}}(d) = (s-1)d - \binom{s+2}{3} - \binom{s+2}{4} + 1,
\]
{\em and $s \ge 2$ is the unique integer for which $\binom{s+2}{3} \le d <
\binom{s+3}{3}$.  Furthermore, if $g = G_{\mbox{\em min}}(d)$, then $s =
s_0(c)$, the least degree of a hypersurface containing $C$.}

\bigskip
\noindent
{\bf Proof.} The simplest way to see this is to consider the $h$-vector of the
curve \cite[\S $1.4$]{M}.  This is a sequence of positive integers $c_0 = 1$,
$c_1,c_2,\dots,c_r$, which determine the degree and genus of the curve
according to the formulae
\[
d = \sum_{i=0}^r c_i,\hskip 0.5 in g = \sum_{i=2}^r (i-1)c_i.
\]
The $c_i$ measure the Hilbert function of a graded ring $R = k[x_0,x_1,x_2]/J$
of finite length, since $C$ has codimension $3$.  Since $R$ is a quotient of a
polynomial ring in three variables, we have $c_i \le \binom{i+2}{2}$ for $i \ge
1$.  The hypothesis $C$ nondegenerate implies $c_1 = 3$.  Thus $d \ge 4$.  For
a given value of $d$,, the genus will be minimized by making each $c_i$ as
large as possible for $i = 2,3,\dots$.  Thus for $4 \le d < 10$ the minimum $g$
is attained by the $h$-vector $1,3,d-4$, with genus $g = d-4$.  For $10 \le d <
20$, the minimum genus comes from the $h$-vector $1,3,6,d-10$, with $g = 2d -
14$.  For the general case, a short calculation with binomial coefficients
gives the formula above.

The least degree $s_0(c)$ of a hypersurface containing $C$ can be read from the
$h$-vector as the least $i$ for which $c_i < \binom{i+2}{2}$.  For the
$h$-vectors giving the minimum genus, this is just the number $s$ defined above.

\bigskip
\noindent
{\bf Remark 3.2.}  It seems reasonable to expect that for each $d \ge 4$ and $g
= G_{\mbox{min}}(d)$, the set of ACM curves of degree $d$ and genus $g$ in
${\mathbb P}^4$ should form an open subset of an irreducible component of the
Hilbert scheme of curves in ${\mathbb P}^4$, and that a general such curve
should be nonsingular, but we do not know how to prove this.

\bigskip
\noindent
{\bf Notation 3.3.}  We will be dealing with curves on certain rational ACM
surfaces in ${\mathbb P}^4$, so here we fix some terminology and notation.

The {\em smooth cubic scroll} $S$ is obtained by blowing up one point $P \in
{\mathbb P}^2$, and embedding in ${\mathbb P}^4$ by the complete linear system
$H = 2l-e$, where $l$ is the total transform of a line in ${\mathbb P}^2$, and
$e$ is the class of the exceptional divisor $E$.  One knows that $\mbox{Pic } S
= {\mathbb Z} \oplus {\mathbb Z}$, generated by $l,e$.  We denote the divisor
class $al-be$ by $(a;b)$.

The {\em Del Pezzo surface} $S$ is obtained by blowing up five points
$P_1,\dots,P_5$, no three collinear in ${\mathbb P}^2$, and embedding in
${\mathbb P}^4$ by $H = 3l - \Sigma e_i$.  In this case $\mbox{Pic } S =
{\mathbb Z}^6$, and we denote the divisor class $al - \Sigma b_ie_i$ by
$(a;b_1,\dots,b_5)$.  If some $b$'s are repeated, we denote that with an
exponent.  Thus in the discussion of $(8,4)$ curves below, the divisor class
$(5;2^2,1^3)$ means $(5;2,2,1,1,1)$.

A {\em Castelnuovo surface} $S$ in ${\mathbb P}^4$ is a smooth surface of
degree $5$ and sectional genus $2$.  It can be obtained by blowing up $8$
points $P_0,P_1,\dots,P_7$ in ${\mathbb P}^2$ and embedding by the linear
system $H = (4;2,1^7)$ (see \cite{Okonek}).  If the points $P_i$ are no three
collinear and no $6$ on a conic, we call it a {\em general} Castelnuovo
surface.  Here $\mbox{Pic } S = {\mathbb Z}^9$, and we denote the divisor class
$al - \Sigma b_ie_i$ by $(a;b_i)$.

A {\em Bordiga surface} is a smooth surface of degree $6$ and sectional genus
$3$.  It can be obtained by blowing up ten points $P_1,\dots,P_{10}$ in
${\mathbb P}^2$ and embedding by the linear system $H = (4;1^{10})$
\cite{Okonek}.  If the points $P_i$ are such that no three are collinear, no $6$
on a conic, and no $10$ on a cubic curve, we call it a {\em general} Bordiga
surface.

\bigskip
\noindent
{\bf Proposition 3.4.}  a) {\em If $C$ is an integral nondegenerate {\em ACM}
curve of degree $d \le 9$ in ${\mathbb P}^4$, the degree-genus pair $(d,g)$
must be one of the following:  $(4,0)$, $(5,1)$, $(6,2)$, $(7,3)$, $(8,4)$,
$(8,5)$,
$(9,5)$, $(9,6)$, $(9,7)$.}

b) {\em For each $(d,g)$ pair in} a), {\em the set of nonsingular nondegenerate
curves in ${\mathbb P}^4$ forms an open subset of an irreducible component of
the Hilbert scheme, and}

c) {\em For each $(d,g)$ pair as above, the general such curve is {\em ACM} and
can be obtained by ascending Gorenstein biliaisons from a line.}

\bigskip
\noindent
{\bf Proof.}  a) A lower bound on $g$ is given by $(3.1)$; an upper bound is
given by the Castelnuovo bound for the genus of an integral curve (see, e.g.,
Rathmann \cite{Rath}).  This list gives all possible values of $g$ between the
lower and upper bounds.

b) For $g = d-4$, the irreducibility is given by a theorem of Ein \cite{Ein}. 
For $(d,g) = (8,5)$, $C$ is the canonical embedding of a non-hyperelliptic
curve of genus $5$, so the family is irreducible.  For $(d,g) = (9,6)$ th curve
is a non-trigonal curve of genus $6$, embedded by a linear system $D = K-P$, so
the family is irreducible (I am indebted to E.~Drozd for this observation). 
For $(d,g) = (9,7)$, the family is irreducible by a theorem of Harris
\cite{EH}.

c) We do a case-by-case analysis.

For $(d,g) = (4,0)$, up to automorphisms of ${\mathbb P}^4$, there is just one
rational normal curve $C$ of degree $4$.  It lies on a smooth cubic scroll $S$,
having divisor class $(2;0)$.  If $H$ denotes the hyperplane class $(2;1)$,
then $C-H = (0;-1)$, which is a line.  Thus $C$ is obtained by an ascending
Gorenstein biliaison from a line on the surface $S$.

For $(d,g) = (5,1)$, suppose given a smooth nondegenerate $(5,1)$ curve $C$ in
${\mathbb P}^4$.  Let $C_0$ be the abstract elliptic curve, and let $D_0$ be
the divisor corresponding to ${\mathcal O}_C(1)$.  Then $C$ is obtained by
embedding $C_0$ with the complete linear system $|D_0|$.  Choose $F$ a divisor
of degree $3$ on $C_0$, and use $|F|$ to embed $C_0$ as a nonsingular cubic
curve $C_1$ in ${\mathbb P}^2$.  Choose a point $P \in C_1$.  Blow up $P$ in
${\mathbb P}^2$ and embed by the linear system $H = 2l - e$ to get a
nonsingular cubic scroll $S$ in ${\mathbb P}^4$.  The image of $C_1$ will be a
$(5,1)$ curve $C_2 \subseteq {\mathbb P}^4$, obtained by embedding $C_0$ with
the linear system $|2F-P|$.  By adjusting the choice of $P$, we may arrange
that $D_0 \sim 2F-P$.  Then $C$ and $C_2$ will differ by an automorphism of
${\mathbb P}^4$.  We conclude that $C$ lies on a smooth cubic scroll $S'$, and
has divisor class $(3;1)$ on $S'$.  Then $C-H = (1;0)$ is a conic.  The conic
in turn can be obtained ascending Gorenstein biliaison from a line on a
plane.  Thus $C$ is obtained by ascending Gorenstein biliaisons from a line.

A similar argument shows that every smooth nondegenerate $(6,2)$ curve $C$ in
${\mathbb P}^4$ appears as a divisor of type $(4;2)$ on a smooth cubic scroll. 
Then $C-H = (2;1)$ is a twisted cubic curve, which can in turn be obtained by
an ascending Gorenstein biliaison on a quadric surface in ${\mathbb P}^3$.

The case $(d,g) = (7,3)$ is a little more complicated.  We will consider only a
smooth nondegenerate non-hyperelliptic $(7,3)$ curve $C$ in ${\mathbb P}^4$. 
Let $C_0$ be the abstract curve of genus $3$, and $D_0$ the divisor giving the
embedding $C$.  Let $C_1$ be the canonical embedding of the non-hyperelliptic
curve $C_0$ as a smooth plane quartic curve in ${\mathbb P}^2$.  Choose five
points $P_1,\dots,P_5$ on $C_1$ with no $3$ collinear.  Blow up $P_1,\dots,P_5$
and embed by $3l - \Sigma e_i$ to get a Del Pezzo surface $S$ in ${\mathbb
P}^4$.  The image of $C_1$ is then a smooth $(7,3)$ curve $C_2 \subseteq S$
with divisor class $(4;1^5)$, which is an embedding of $C_0$ by the divisor $3K
- \Sigma P_i$.  We would like to choose the $P_i$ so that $3K - \Sigma P_i \sim
D_0$, i.e., $\Sigma P_i \sim 3K - D_0$.

{\bf Case 1.}  The divisor $3K - D_0$ can be represented by $5$ points, no
three collinear.  In this case we find that $C$ is contained in a Del Pezzo
surface $S'$, with divisor class $(4;1^5)$.  Then $C-H = (1;0^5)$, which is a
twisted cubic curve, so $C$ can be obtained from a line by ascending Gorenstein
biliaisons.

{\bf Case 2.}  If Case~1 does not occur, one sees easily that $3K-D_0 \sim K +
P$ for some point $P$.  In this case we take $S$ to be the smooth cubic scroll
obtained by blowing up $P$.  Then $C_2$ is a divisor of type $(4;1)$ on $S$,
and is an embedding of $C_0$ by $2K - P \sim D_0$.  So $C$ lies on a cubic
scroll, and $C - 2H = (0;-1)$ is a line.

Note the two types of non-hyperelliptic $(7,3)$ curves can be distinguished by
the property that the Case~1 curves have only finitely many trisecants, while
the Case~2 curves have infinitely many trisecants.

For a smooth $(8,4)$ curve $C$ we use a different technique.  For this curve,
$h^0({\mathcal I}_C(2)) = 2$, so $C$ is contained in a unique complete
intersection surface $S = F_2 \cap F'_2$.  Since the family of $(8,4)$ curves
is irreducible, and since a general complete intersection surface $S = F_2 \cap
F'_2$ is a smooth Del Pezzo surface containing a smooth $(8,4)$ curve $C$ in
the divisor class $(5;2^2,1^3)$, we conclude that a general such $C$ lies on a
smooth Del Pezzo surface $S$, with divisor class $(5;2^2,1^3)$.  Then $C-H =
(2;1^2,0^3)$, which is a nondegenerate smooth $(4,0)$ curve in ${\mathbb
P}^4$.  Thus using the case of $(4,0)$ above, we conclude that $C$ can be
obtained by ascending Gorenstein biliaisons from a line.

The case $(9,6)$ is similar to $(8,4)$, because again $h^0({\mathcal I}_C(2)) =
2$, we conclude that a general smooth $(9,6)$ curve $C$ lies on a Del Pezzo
surface $S$ with divisor class $(6;2^4,1)$.  In this case $C-2H = (0;0^4,-1)$
which is a line.

A smooth nondegenerate $(8,5)$ curve $C$ in ${\mathbb P}^4$ is the canonical
embedding of a non-hyperelliptic genus $5$ curve.  According to the theorem of
Petri \cite{S-D} if the curve is not trigonal, then $C$ is the complete
intersection of three quadric hypersurfaces $C = F_2 \cap F'_2 \cap F''_2$. 
Let $S = F_2 \cap F'_2$.  Then $C$ is the divisor $2H$ on $S$, and $C-H$ is an
elliptic quartic curve in ${\mathbb P}^3$, which can be obtained from a conic
by biliaison on a quadric surface.  Thus $C$ is obtained by ascending
Gorenstein biliaisons from a line.

If the curve is trigonal, then $C$ lies on a smooth cubic scroll $S = F_2 \cap
F'_2 \cap F''_2$.  It has divisor class $(5;2)$, so $C - 2H$ is $(1;0)$, a
conic, and we conclude again.  Note in this case $C = H - K$ is arithmetically
Gorenstein, even though it is not a complete intersection.

For $(d,g) = (9,5)$, arguments like the ones above show that we can embed the
general genus $5$ curve as a plane quintic with a double point, and thus obtain
a general $(9,5)$ curve $C$ on a Castelnuovo surface with divisor class
$(5;2,1^7)$.  Then $C - H = (1,0^8)$ is a $(4,0)$ curve and we use our earlier
result.

Finally, every smooth $(9,7)$ curve $C$, as a curve of maximal genus, lies on a
cubic surface $S$ as the divisor $3H$,  by a theorem of Harris \cite{EH}. 
Then $C-2H$ is a twisted cubic curve, and we are done.

\bigskip
\noindent
{\bf Corollary 3.5.}  {\em A general smooth {\em ACM} curve of degree $d \le 9$
in
${\mathbb P}^4$ is glicci.}

\bigskip
\noindent
{\bf Remark 3.6.}  The glicciness of ACM curves lying on general ACM surfaces
in ${\mathbb P}^4$ and of integral ACM curves of degree $\le 7$ was already
proven in \cite[\S 8]{KMMNP}.  Our contribution is to show that a general
smooth ACM curve of degree $\le 9$ actually does lie on a smooth rational ACM
surface, and to check the possibly stronger property that they can be obtained
from a line by ascending Gorenstein biliaisons (cf.\ Question $1.3$b).

Note that our proof actually shows {\em every} smooth curve with $(d,g) =
(4,0),(5,1),(6,2)$ is ACM and can be obtained by ascending Gorenstein
biliaisons from a line.  The same can be said for $(7,3)$ curves, by extending
the analysis above:  one can show that a hyperelliptic $(7,3)$ cuve lies on a
cubic scroll or the cone over a twisted cubic curve, and in both cases is
obtained from a smooth $(4,0)$ curve by biliaison.

For the next case of $(8,4)$ curves, the situation is more complicated.  There
are smooth hyperelliptic $(8,4)$ curves on a cubic scroll, but they are not
ACM.  Since an ACM $(8,4)$ curve lies on a unique complete intersection surface
$S = F_2 \cap F'_2$, to study  {\em all} smooth $(8,4)$ curves, one would
presumably have to study the possible singular surfaces $S$.  One approach is
to use Riemann--Roch on the surface to show that the divisor $C-H$ is
effective, but then one has to deal with not necessarily irreducible $(4,0)$
curves.

The analysis becomes increasingly complex for the remaining cases, so we do not
know if c) holds for all smooth ACM curves of the given degree and genus.

\bigskip
\noindent
{\bf Remark 3.7.}  For $d \ge 10$, the family of smooth non-degenerate ACM
curves of given $(d,g)$ in ${\mathbb P}^4$ may not be irreducible.  The first
example is $(d,g) = (10,9)$, for which there are two different families of such
curves lying on smooth cubic scrolls.

\bigskip
\noindent
{\bf Example 3.8.}  We consider smooth $(10,6)$ curves in ${\mathbb P}^4$.  Note
that all the curves in Proposition $3.4$, being ACM of degree $\le 9$, are
contained in quadric hypersurfaces, since their hyperplane section is $\le 9$
points and is contained in a quadric surface of ${\mathbb P}^3$.  The case
$(d,g) = (10,6)$ is the first case where there are smooth ACM curves not
contained in a quadric hypersurface.

By the theorem of Ein \cite{Ein}, the family of smooth $(10,6)$ curves in
${\mathbb P}^4$ is irreducible.  To show that a general $(10,6)$ curve is ACM,
it suffices, by semicontinuity, to exhibit one.  A divisor of type
$(5;1^{10})$ on a general Bordiga surface $S$ is the transform of a plane
quintic curve, which can be taken to be smooth, so we get a smooth $(10,6)$
curve $C$ on $S$.  For this curve $C-H = (1;0^{10})$ is a smooth $(4,0)$
curve, which is ACM, so $C$ is also ACM.  Note however that the curve just
described is not general in the variety of moduli of curves of genus $6$,
because it has a $g_5^2$:  a representation as a plane quintic curve.

Next, let $C_0$ be an abstract curve of genus $6$, with general moduli.  Then
$C_0$ admits a birational representation as a plane curve $C_1 \subseteq
{\mathbb P}^2$ with four nodes $P_1,P_2,P_3,P_4$, no three collinear
\cite{ACGH}.  Choose six additional points $P_5,\dots,P_{10}$ on $C_1$ in
general position.  Blow up $P_1,\dots,P_{10}$ to obtain a Bordiga surface $S$,
containing the proper transform $C_2 \subseteq S$ of $C_1$.  Then $C_2$ is a
smooth $(10,6)$ curve in ${\mathbb P}^4$ with general moduli.  Since the curve
has genus $6$, by varying the choice of the six points $P_5,\dots,P_{10}$, we
can obtain any general divisor class on $C_2$ as its hyperplane section.  We
conclude that the general $(10,6)$ curve $C$ in ${\mathbb P}^4$ is contained
in a general Bordiga surface $S$ with divisor class $(6;2^4,1^6)$.  Then $C-H
= (2;1^4,0^6)$, which is a smooth $(4,0)$ curve, so $C$ can be obtained by
ascending Gorenstein biliaisons from a line.

To study the $(10,6)$ curves in more detail, we note that as a general Bordiga
surface $S$, there are eight divisor classes (up to permutation of the $P_i$)
containing $(10,6)$ curves.  They are
\begin{eqnarray*}
D_1 &= &(5;1^{10}) \\
D_2 &= &(6;2^4,1^6) \\
D_3 &= &(7;2^9,0) \\
D_4 &= &(7;3,2^6,1^3) \\
D_5 &= &(8;3^3,2^6,1) \\
D_6 &= &(8;4,2^9) \\
D_7 &= &(9;3^6,2^4) \\
D_8 &= &(10;3^{10}).
\end{eqnarray*}
Of these $D_3$ and $D_6$ have Rao module $k$.  They will be discussed in the
next section.  The remaining $6$ cases are ACM.  The first three of these,
$D_1$, $D_2$, and $D_4$, can be obtained by Gorenstein biliaison from $(4,0)$
curves on $S$.  However, $D_5 - H$, $D_7 - H$, and $D_8 - H$ are not
effective divisors so these curves cannot be obtained by Gorenstein biliaison
on this surface $S$.

Using the arithmetically Gorenstein divisor $3H - K$ on $S$, of degree $20$,
the divisor class $D_i$ is Gorenstein-linked to $D_{8-i}$.  It follows that
$D_5$, $D_7$, $D_8$ are glicci (as observed in \cite[\S 8]{KMMNP}).  However, we
do not know whether or not these curves may be obtained by ascending Gorenstein
biliaison on some other surface.

\bigskip
\noindent
{\bf Example 3.9.}  For our last example, we will study ACM $(20,26)$ curves
in ${\mathbb P}^4$.  Note that all the curves in the earlier part of this
section, plus all the ACM curves lying on rational ACM surfaces in ${\mathbb
P}^4$, which were proved to be glicci in \cite[\S 8]{KMMNP}, lie on cubic
hypersurfaces in ${\mathbb P}^4$.  So by analogy with our findings for points
in ${\mathbb P}^3$ in \S 2 above, we might expect that all ACM curves contained
in cubic hypersurfaces in ${\mathbb P}^4$ would be glicci.  This also suggests
that in looking for counterexamples to ACM $\Rightarrow$ glicci (Question
$1.3$), we should look at curves not contained in a cubic hypersurface.

The first example of an ACM curve in ${\mathbb P}^4$ not contained in a cubic
hypersurface will have $h$-vector $1,3,6,10$ (cf.\ proof of $3.1$).  It has
degree $20$ and genus $26$.  For existence of such curves, we let $C$ be the
determinantal curve defined by the $4 \times 4$ minors of a $4 \times 6$ matrix
of general linear forms.  A general such curve will be smooth, ACM, of degree
$20$ and genus $26$.  The family of such determinantal curves has dimension
$\le 69$, by \cite[10.3]{KMMNP}.

The method of \cite[3.7]{KMNP} shows that $C$ is linearly equivalent to $H+K$
on an ACM surface $S$ in ${\mathbb P}^4$, of degree $10$ and sectional genus
$11$, defined by the $4 \times 4$ minors of a $4 \times 5$ matrix of general
linear forms, where $H$ denotes the hyperplane section of $S$, and $K$ denotes
the canonical class of $S$.  Furthermore, a similar argument using
\cite[3.1]{KMNP} shows that the curve $C_0$ defined by the $3 \times 3$ minors
of a $3 \times 5$ matrix of linear forms will be linearly equivalent to $K$ on
$S$.  This latter curve $C_0$ also appears in the divisor class $2H_0+K_0$ on
the surface $S_0$ defined by $3 \times 3$ minors of a $3 \times 4$ matrix of
linear forms.  (I am grateful to J.~Migliore for pointing out the paper
\cite{KMNP} and explaining to me how to obtain these linear equivalences.)

Now $S_0$ is just the Bordiga surface, and $C_0$ is an ACM $(10,6)$ curve,
discussed earlier.  By the linear equivalence $C_0 \sim 2H_0 + K_0$ on $S_0$ we
recognize that $C_0$ is in the class $(5;1^{10})$, which we called $D_1$ in
$(3.8)$ above.  These are isomorphic to plane curves of degree $5$ and thus are
not general in the moduli of curves of genus $6$.

Since $C_0$ can be obtained by ascending Gorenstein biliaison from a line, and
since $C \sim C_0 + H = K + H$ as the ACM surface $S$, we conclude that the
determinantal $(20,26)$ curve $C$ can also be obtained by ascending Gorenstein
biliaison from a line.

Next, I claim that the only way to obtain an ACM $(20,26)$ curve $D$ in
${\mathbb P}^4$ by ascending Gorenstein biliaison is from an ACM $(10,6)$ curve
$C_1$ on an ACM surface $S_1$ of degree $10$ and sectional genus $11$, as $D
\sim C_1 + H$ on $S_1$.  Indeed, suppose that $D \sim C_1 + H$ on some ACM
surface $S_1$.  Then $C_1$ is an ACM curve of type $(d_1,g_1)$ in ${\mathbb
P}^4$, while $H$ is an ACM curve of type $(d_2,g_2)$ in ${\mathbb P}^3$.  From
this we get $(20,26) = (d_1 + d_2, g_1 + g_2 + d_1 -1)$.  For each $d_1$
(resp.~$d_2$) we know the minimum possible genus of an ACM curve in ${\mathbb
P}^4$ (resp.~${\mathbb P}^3$)---cf.\ $3.1$.  Looking at these, we find that
$g_1 + g_2 + d_1 - 1 > 26$ in all cases except $(d_1,g_1) = (10,6)$ and
$(d_2,g_2) = (10,11)$.  Thus any $(20,26)$ curve $D$ that can be obtained by
ascending Gorenstein biliaison must lie on a surface $S_1$ of degree $10$ and
sectional genus $11$.

Now we look at the dimensions of some families of $(20,26)$ curves.  By
\cite[$10.3$]{KMMNP}, the family of determinantal curves $C$ as above has
dimension $\le 69$.  On the other hand, each component of the Hilbert scheme of
$(20,26)$ curves in ${\mathbb P}^4$ has dimension $\ge 5d + 1-g = 75$.  So we
see immediately that a general element of an irreducible component of
$H_{20,26}$ cannot be determinantal.  However, there may be other $(20,26)$
curves $C'$ on $S$, not determinantal themselves, but linearly equivalent to
$C$, obtainable by ascending Gorenstein biliaison on $S$.

So let us find the dimension of the complete linear system $|C|$ on $S$.  From
the exact sequence
\[
{\mathcal O} \rightarrow {\mathcal O}_S \rightarrow {\mathcal O}_S(C) 
\rightarrow {\mathcal O}_C(C)
\rightarrow 0
\]
we see that $\dim_S|C| = h^0({\mathcal O}_C(C)) = C^2 + 1-g + h^1({\mathcal
O}_C(C))$.  We also have a resolution of ${\mathcal O}_S$
\[
0 \rightarrow {\mathcal O}_{{\mathbb P}^4}(-5)^4 \rightarrow {\mathcal
O}_{{\mathbb P}^4} (-4)^5 \rightarrow {\mathcal O}_S \rightarrow 0
\]
coming from its matrix representation.  From this we find $h^2({\mathcal O}_S)
= 4$ and $p_a(S) = 4$.  On the surface $S$ we have $H^2 = \deg S = 10$.  From
the adjunction formula for $H$, which is a $(10,11)$ curve, we find $H\cdot K =
10$.  And from the formula \cite[p.~434]{AG} for surfaces in ${\mathbb
P}^4$, we find $K^2 = 5$.  Now $C = H + K$, so we get $C^2 = 35$.  Also, since
$C = H + K$, from the Kodaira Vanishing Theorem we have $H^1({\mathcal O}_S(C))
= H^2({\mathcal O}_S(C)) = 0$.  Thus $h^1({\mathcal O}_C(C)) \cong
h^2({\mathcal O}_S) = 4$.  So we find
\[
\dim_S|C| = 35 + 1 - 26 + 4 = 14.
\]
The family of ACM  surface $S$ has dimension $60$ (for example by the formula
of Ellingsrud \cite{E}), so we find that the family of $(20,26)$ curves in
${\mathbb P}^4$ that can be linearly equivalent to $C$ on such a surface $S$
has dimension $\le 74$.  In particular, a general curve in an irreducible
component of $H_{20,26}$ does not arise in this way.

There may also be other linear equivalence classes of $(20,26)$ curves on $S$
of larger dimension.  Indeed, this is what does happen with the $(10,6)$ curves
studied in $(3.8)$ above:  the curves linearly equivalent to the determinantal
curves were all of type $D_1$ on the Bordiga surface, and these curves are not
general in the moduli of genus $6$ curves, while a general genus $6$ curve
appears as an ACM curve in a different linear system $D_2$ on the Bordiga
surface.  So we must see if something analogous happens with the $(20,26)$
curves.

First, we look on a general ACM surface $S$ of degree $10$ and sectional genus
$11$.  According to a theorem of Lopez \cite[III.4.2]{Lop}, $\mbox{Pic } S =
{\mathbb Z} \oplus {\mathbb Z}$ generated by $H$ and $K$.  We look for divisors
$mH + nK$ with degree $20$ and genus $26$.  There are only two
possibilities:  $C = H + K$ or $C' = 4H - 2K$.  In the latter case we compute
$C'{}^2 = 20$.  Therefore, by Clifford's theorem, $h^0({\mathcal O}_{C'}(C')) -
1 \le 10$, so $\dim_S|C'| \le 11$.  Thus the family of such curves $C'$ in
${\mathbb P}^4$ has dimension $\le 71$.  So we see that a general ACM $(20,26)$
curve in ${\mathbb P}^4$ cannot lie on a {\em general} ACM surface $S$ of
degree $10$ and sectional genus $11$.

Now let us estimate the dimension of a family of smooth $(20,26)$ curves $D$,
general in an irreducible component of $H_{20,26}$ containing the determinantal
curves $C$ above, and lying on a non-general ACM surface $X$ of degree $10$ and
sectional genus $11$.  We will make use of the Clifford index of a curve.

Recall that the {\em gonality} of a curve $C$ is the least $d$ for which there
exists a linear system $g^1_d$ on the curve.  The {\em Clifford index} of the
curve is the minimum of $d-2r$, taken over all linear systems $g^r_d$ with $r
\ge 1$ and $0 < d \le g-1$.  For most curves, the Clifford index is equal to
$\mbox{gon}(C) - 2$, computed by a $g^1_d$.  Curves for which this is not so are
{\em Clifford exceptional} curves, and have been studied by Martens \cite{Mar}
and Eisenbud et al.\ \cite{ELMS}.

If $C$ is the determinantal $(20,26)$ curve studied above, then $C \sim H + K$
on the surface $S$.  The hyperplane section $H$ is a $(10,11)$ curve in
${\mathbb P}^3$, obtained as $C_0 + H_0$ on a nonsingular quartic surface in
${\mathbb P}^3$.  Here $C_0$ is a nonhyperelliptic $(6,3)$ curve having
gonality $3$; $H_0$ is a plane quartic curve, also having gonality $3$.  Hence,
by \cite{CIA}, $H$ has gonality $\ge 6$.  (In fact the gonality
is equal to $6$ because $H$ must have a $4$-secant.)  The curve $K$ on $S$ is
the determinantal $(10,6)$ curve discussed above, isomorphic to a plane quintic
curve, with gonality $4$.  So applying \cite{CIA} again, we
find the gonality of $C$ is $\ge 10$.  It follows from the study of Clifford
exceptional curves in \cite{Mar} and \cite{ELMS} that $C$ is not exceptional, so
we conclude $\mbox{Cliff } C \ge 8$.  (On the other hand, the linear system
$|K|$ on $S$ cuts out a $g_{15}^3$ on $C$, so $\mbox{Cliff } C \le 9$.  I
suspect $\mbox{Cliff } C = 9$, but don't know how to prove that.)  Since we are
considering a curve $D$ that is general in an irreducible component of
$H_{20,26}$ containing $C$, we may assume also that $\mbox{Cliff } D \ge 8$.

Now we consider a smooth $(20,26)$ curve $D$ with $\mbox{Cliff } D \ge 8$,
contained in a smooth ACM surface of degree $10$ and sectional genus $11$ in
${\mathbb P}^4$, and we want to estimate the dimension of the linear system
$|D|$ on $S$.  As above, we find
\[
\dim_S|D| = D^2 + 1 - g + h^1({\mathcal O}_D(D)).
\]
Since $D$ is a $(20,26)$ curve, the adjunction formula gives 
\[
D^2 + D \cdot K = 50.
\]
Let us denote $D \cdot K$ by $b$.  Then $D^2 = 50 - b$.  On the other hand, let
us consider the linear system $|D \cdot K|$ on $D$.  It has dimension $a-1$,
where $a = h^0({\mathcal O}_D(K))$.  Since $K_D = (D+K)\cdot D$, we also have
$h^1({\mathcal O}_D(D)) = a$.  The linear system $|D\cdot K|$ thus has
dimension $a-1$ and degree $b$.  Our hypothesis $\mbox{Cliff } D \ge 8$ thus
implies $b - 2a + 2 \ge 8$, or $b \ge 2a + 6$.

Now we can compute
\begin{eqnarray*}
\dim_S|D| &= &D^2 + 1 - g + h^1({\mathcal O}_D(D)) \\
&= &50 - b + 1 - 26 + a \\
&= &25 + a - b.
\end{eqnarray*}
Since $b \ge 2a + 6$, we find
\[
\dim_S|D| \le 19 - a.
\]
Now from the exact sequence
\[
0 \rightarrow {\mathcal O}_X(K-D) \rightarrow {\mathcal O}_X(K) \rightarrow
{\mathcal O}_D(K) \rightarrow 0,
\]
we find $a = h^0({\mathcal O}_D(K)) \ge h^0({\mathcal O}_X(K)) = h^2({\mathcal
O}_X) = 4$.  Thus
\[
\dim_S|D| \le 15.
\]
On the other hand, our surface $S$ is not general, so it moves in a family of
dimension at most $59$, so the dimension of the family of curves that arise in
this way is at most $74$.

In conclusion, we see that there exists an irreducible component of the Hilbert
scheme $H_{20,26}$ of smooth $(20,26)$ curves in ${\mathbb P}^4$ (namely one
containing the determinantal curves) whose general member is an ACM curve that
does not lie on an ACM surface $S$ of degree $10$ and sectional genus $11$, and
so cannot be obtained by ascending Gorenstein biliaison from a line.  We
propose this curve as a possible candidate for a counterexample to ACM
$\Rightarrow$ glicci (Question $1.3$).

\section{Curves in ${\mathbb P}^4$ with Rao module $k$}

Let ${\mathcal M}$ be the set of all locally CM curves in ${\mathbb P}^4$ with
Rao module $k$ (i.e., of dimension one in one degree only).  One knows that the
Rao module must occur in a nonnegative degree \cite[$1.3.11$(b)]{M}, and that
there are curves with Rao module $k$ in degree $0$ (e.g., two skew lines).  So
we denote by ${\mathcal M}_h$ the set of curves with Rao module $k$ in degree
$h$, and note that ${\mathcal M} = \cup_{h \ge 0} {\mathcal M}_h$.

Let ${\mathcal L} \subseteq {\mathcal M}$ be the subset of those curves in the
$G$-liaison class of two skew lines, and let ${\mathcal L}_h = {\mathcal L}
\cap {\mathcal M}_h$.  Then ${\mathcal L} = \cup_{h \ge 0} {\mathcal L}_h$, and
the curves in ${\mathcal L}_0$ are the {\em minimal} curves defined in \S1
above.

In this section we will study the curves in ${\mathcal M}$, with a view to
elucidating Questions $1.4$ and $1.5$ above.

\bigskip
\noindent
{\bf Proposition 4.1.} a) {\em ${\mathcal M}_0$ contains curves of every degree
$d \ge 2$.}

b) {\em For each $d \ge 2$, the set of curves in ${\mathcal M}_0$ of degree
$d$ forms an irreducible family, whose general member is the disjoint union $C
= C' \cup L$ of a plane curve $C'$ of degree $d-1$ and a line $L$, not meeting
the plane of $C'$.}

c) {\em Every curve in ${\mathcal M}_0$ is in the $G$-liaison class of two
skew lines, i.e.\ ${\mathcal L}_0 = {\mathcal M}_0$.}

\bigskip
\noindent
{\bf Proof.} a) The case of two skew lines in ${\mathbb P}^3$ is well-known
\cite[Example $6.2$, p.~34]{MDP}.  For $d \ge 3$, let $C = C' \cup L$ as
described in b).  Clearly $h^0({\mathcal O}_C) = 2$, so $h^1({\mathcal I}_C) =
1$.  On the other hand, since $C'$ and $L$ are contained in disjoint sublinear
spaces of
${\mathbb P}^4$, it is clear that $H^0({\mathcal O}_{{\mathbb P}^4}(n)) \rightarrow
H^0({\mathcal O}_C(n))$ is surjective for $n \ge 1$, so $C \in {\mathcal M}_0$.

b) I claim any degree~$2$ curve $C$ in ${\mathbb P}^4$ with $M = k$ lies in
${\mathbb P}^3$.  If the curve is reduced, it is two lines, hence in a
${\mathbb P}^3$.  If it is not reduced, then it is a double structure on a line
$L$, and we have an exact sequence
\[
0 \rightarrow {\mathcal L} \rightarrow {\mathcal O}_C \rightarrow {\mathcal
O}_L \rightarrow 0
\]
where ${\mathcal L}$ is an invertible sheaf on $L$.  Then there is a surjective
map $u: {\mathcal I}_L/{\mathcal I}_L^2 \rightarrow {\mathcal L} \rightarrow
0$, and ${\mathcal L} \cong {\mathcal O}_L(a)$ for some $a$.  Since ${\mathcal
I}_L/{\mathcal I}_L^2 \cong {\mathcal O}_L(-1)^3$, the map $u$ is given by
three sections of ${\mathcal O}_L(a+1)$.  If $a = -1$, we get a double line in
a plane.  If $a=0$, there is a linear form $x$ killed by $u$, so $C$ lies in
the ${\mathbb P}^3$ defined by $x=0$.  If $a > 0$, then the exact sequence
\[
H^0({\mathcal O}_L(-1))^3 \rightarrow H^0({\mathcal O}_L(a)) \rightarrow
H^1({\mathcal I}_C) \rightarrow 0
\]
shows the Rao module is bigger than $k$.

Thus a curve of degree~$2$ with $M=k$ lies in a ${\mathbb P}^3$, so these form
an irreducible family whose general member is two skew lines.

So now let $d \ge 3$.  Then $C$ cannot be contained in ${\mathbb P}^3$, because
of the Lazarsfeld--Rao property for curves in ${\mathbb P}^3$, so
$h^0({\mathcal O}_C) = 2$, because of the Rao module, and $h^0({\mathcal
O}_C(1)) = 5$.  Let $A = H_*^0({\mathcal O}_C)$.  This is a graded $S$-algebra,
where $S = k[x_0,x_1,x_2,x_3,x_4]$, and in particular $A_0$ is a
$2$-dimensional $k$-algebra.  We consider two cases.

{\bf Case 1.}  $A_0$ is reduced, hence isomorphic to $k \times k$ as a
$k$-algebra.  Then $A_0$ contains two orthogonal idempotents $e',e''$, such
that $e' + e'' = 1$, $e'{}^2 = e'$, $e''{}^2 = e''$, and $e'e'' = 0$.  Hence
$C$ is the disjoint union of two curves $C',C''$, defined by the vanishing of
$e',e''$, respectively.  Let $H',H''$ be the linear spans of the curves
$C',C''$.  Then $h^0({\mathcal O}_C(1)) = h^0({\mathcal O}_{H'}(1)) +
h^0({\mathcal O}_{H''}(1)) = 5$.  So one of these, say $H'$, is a plane, and
the other, $H''$ is a line $L$.  Thus $C'$ is a plane curve in $H'$, and $C =
C' \cup L$ as required.  Note that $H',H''$ do not meet since $h^0({\mathcal
O}_{{\mathbb P}^4}(1)) = h^0({\mathcal O}_{H'}(1)) + h^0({\mathcal
O}_{H''}(1))$.

{\bf Case 2.}  $A_0$ is non-reduced, in which case it is isomorphic to the ring
$k[\epsilon]/(\epsilon^2)$.  Let $f \in A_0$ be a nonzero element with $f^2 =
0$.  Now $A_1 \cong S_1$ is the $k$-vector space generated by
$x_0,x_1,x_2,x_3,x_4$.  Multiplication by $f$ on $A_1$ is a nilpotent linear
map with $f^2 = 0$.  Furthermore, since $C$ is locally CM, the kernel of $f$
acting on $A_1$ must have dimension $\le 3$.  Otherwise $f$ would be supported
at a point.  So $f$ has rank $\ge 2$.  Now from the structure of nilpotent
transformations it follows (after a linear change of coordinates) that $fx_0 =
x_2$, $fx_1 = x_3$, $fx_2 = fx_3 = fx_4 = 0$.  Hence we can identify the
$S$-algebra $A$ as
\[
A \cong S[f]/((f^2,fx_0-x_2,fx_1-x_3,fx_2,fx_3,fx_4) + I_C) 
\]
where $I_C \subseteq S$ is the homogeneous ideal of $C$.

Now let $H'$ be the plane $x_2 = x_3 = 0$, and let $C'$ be the curve obtained
from $C \cap H$ by removing its embedded points, if any.  Then there is an
exact sequence
\[
0 \rightarrow {\mathcal L} \rightarrow {\mathcal O}_C \rightarrow {\mathcal
O}_{C'} \rightarrow 0.
\]
Since $C'$ is a plane curve, $h^0({\mathcal O}_{C'}) = 1$, and so
$h^0({\mathcal L}) = 1$.  Furthermore note that the image of $f$ in ${\mathcal
O}_{C'}$ is annihilated by $x_0,x_1,x_2,x_3,x_4$, hence is $0$.  So $f$
generates $h^0({\mathcal L})$.  Now $f$ is annihilated by $x_2,x_3,x_4$, so it
has support on the line $L: x_2 = x_3 = x_4 = 0$.  Thus ${\mathcal L}$ is an
${\mathcal O}_L$-module, it is torsion-free since $C$ is locally CM, and
contains the submodule ${\mathcal O}_L$ generated by $f$.  Hence ${\mathcal L}
\cong {\mathcal O}_L$, generated by $f$.

Now it is clear that $C$ consists of the plane curve $C'$ of degree $d-1$,
containing the line $L$, plus a multiplicity two structure on $L$ with $p_a =
-1$.  This is the limit of a flat deformation of the disjoint unions $C' \cup
L$ described above, as the skew line $L$ approaches a line in the curve $C'$.

So the curves in ${\mathcal M}_0$ of any degree $d \ge 2$ form an irreducible
family.

c) Let $C \in {\mathcal M}_0$ have degree $d$.  The case $d = 2$ in ${\mathbb
P}^3$ is well-known, so we may assume $d \ge 3$.  First consider the disjoint
union $C = C' \cup L$ as in b).  Take a hyperplane ${\mathbb P}^3$ containing
$L$ and meeting the plane $H'$ of $C'$ in a line $L'$, skew to $L$, and not a
component of $C'$.  Let
$Q$ be a nonsingular quadric surface in that
${\mathbb P}^3$ containing $L$ and $L'$.  Then $S = H' \cup Q$ is an ACM
surface of degree $3$ in ${\mathbb P}^4$.  Note that its negative canonical
divisor $-K$ consists of a conic in $H'$ plus a divisor of bidegree $(1,2)$ on
$Q$, meeting
$L'$ in the same two points as the conic, where $(1,0)$ is the class of $L$. 
(Here we leave some details to the reader.)  Now given $C$, there is an AG
divisor
$X$ in the linear system
$(d-3)H-K$ on
$S$ containing $C$ \cite[$4.2.8$]{M}.  The linked curve $D$ is a divisor of
bidegree
$(d-3,d-1)$ on
$Q$, which is in the biliaison class of two skew lines on $Q$.  Thus $C$ is in
${\mathcal L}_0$.

In the special case where $C$ is a plane curve $C'$ containing a line $L$, plus
a double structure on $L$ as above, we use exactly the same construction,
except that now the hyperplane ${\mathbb P}^3$ meets $H'$ in $L$, and the
quadric surface $Q$ contains the double structure on $L$.  The same liaison
works, using the theory of generalized divisors \cite{GD}.

\bigskip
\noindent
{\bf Remark 4.2.} The fact that ${\mathcal L}_0$ is not a single irreducible
family was observed by Migliore \cite[$5.4.8$]{M}, who gave the example of a
curve of degree $3$ in ${\mathcal L}_0$.  His student Lesperance \cite{Les} has
independently proved $3.1$a), c) in the case of reduced curves.  Lesperance has
also shown \cite[$4.5$]{Les} that for other Rao modules, the set of minimal
curves of given degree need not be irreducible.  Thus $(4.1b)$ is special to
the case of Rao module $M=k$.

\bigskip
\noindent
{\bf Example 4.3.}  Let $C$ be a smooth curve of degree $5$ and genus $0$ in
${\mathbb P}^4$, not contained in any ${\mathbb P}^3$.  It is the projection of
the rational normal curve $\Gamma$ in ${\mathbb P}^5$ from a point not lying on
any secant line of $\Gamma$.  A little elementary geometry shows that $C$ has a
unique trisecant $E$.  If $C$ meets $E$ in three distinct points, then the three
points of intersection of
$C$ and
$E$ determine a unique isomorphism (of abstract ${\mathbb P}^1$'s) from $C$ to
$E$ fixing those three points.  Let $S$ be the surface formed as the closure of
the set of lines joining corresponding points of $C$ and $E$.  Then $S$ is a
rational cubic scroll in ${\mathbb P}^4$.

On $S$, our rational quintic $C$ has divisor class $(4;3)$.  The linear system
$C-H = (2;2)$ contains a disjoint union of two rulings of the surface $S$. 
Hence $C$ is obtained by one elementary $G$-biliaison from two skew lines.  In
particular, $C \in {\mathcal L}_1$.

If $E$ is a degenerate trisecant, i.e., a tangent line meeting the curve again,
or an inflectional tangent, we can still show $C \in {\mathcal L}_1$ as
follows.  The smooth $(5,0)$ curves in ${\mathbb P}^4$ form an irreducible
family, so $C$ is a specialization of the general type described above.  Hence
$C$ must lie on a cubic surface in ${\mathbb P}^4$.  It cannot lie on a
reducible surface, since $C$ is not in ${\mathbb P}^3$.  The only other
irreducible cubic surface is the cone over a twisted cubic curve, and that
surface contains no smooth $(5,0)$ curves.  Hence $C$ is on a smooth rational
cubic scroll, and the previous argument applies.

\bigskip
\noindent
{\bf Example 4.4.}  We consider smooth curves of type $(6,1)$ (degree $6$ and
genus $1$) in ${\mathbb P}^4$, not contained in any ${\mathbb P}^3$.  Then
$h^1({\mathcal I}_C(1)) = 1$ and $h^0({\mathcal I}_C(2)) \ge 3$.

{\bf Case 1.}  If three quadric hypersurfaces containing $C$ intersect in a
surface, then that surface must be a cubic rational scroll $S$ (reason:  the
degree of $S$ must be $\le 3$; $C$ is not contained in a plane or a quadric
surface, and there is no $(6,1)$ curve on the cone over a twisted cubic
curve).  In this case $C = (3;0)$ on $S$ and $C-H = (1;-1)$, which contains the
disjoint union of a conic and a line.  Thus $C$ is in ${\mathcal L}_1$ and is
obtained by a single elementary $G$-biliaison from a curve of degree $3$ in
${\mathcal L}_0$.  This curve $C$ has infinitely many trisecants, formed by the
rulings of $S$.

{\bf Case 2.}  Three quadric hypersurfaces containing $C$ will intersect in a
complete intersection curve $X$ of degree $8$ and genus $5$.  The residual
intersection $D$ will be a curve of degree $2$.  $D$ cannot be a plane curve,
because then it would meet $C$ in $5$ points, and projection from the plane of
$D$ would be a birational map of $C$ to a line, which is impossible.  Hence $D$
is two skew lines or a nonplanar double structure on a line.  By reason of the
genus of $X$, $D$ will be either two trisecants of $C$ or a single trisecant. 
Note also that these are all the trisecants of $C$, because any trisecant of
$C$ must be contained in each quadric hypersurface containing $C$, hence in $X$.

{\bf Case 2a.}  An example of a $(6,1)$ curve with two trisecants can be
obtained on a Del Pezzo surface, as a divisor of type $(3;1^3,0^2)$.  In this
case $C - H = (0;0^3,-1^2)$, which is a disjoint union of two lines, so $C$ is
in ${\mathcal L}_1$ and is obtained by one elementary $G$-biliaison from the
minimal curve of degree $2$ in ${\mathcal L}_0$.  This curve $C$ has two
trisecants, the lines $F_{45} = (1;0^3,1^2)$ and $G = (2;1^5)$.

{\bf Case 2b.}  An example of a $(6,1)$ curve with one trisecant can be
obtained as follows.  We project the Veronese surface $V$ in ${\mathbb P}^5$
from a point in a plane containing a conic of $V$, so as to obtain a quartic
surface $S$ in ${\mathbb P}^4$ with a double line $L$.  A general cubic curve
in ${\mathbb P}^2$ gives a $(6,1)$ curve in $V$ meeting the conic in three
points that project to distinct points of the line $L$ in $S$.  Thus the image
$C \subseteq S$ of this curve will be a smooth $(6,1)$ curve having $L$ as a
trisecant.  Now the surface $S$ is smooth except for a double line and two
pinch points, hence locally $CM$.  Its general hyperplane section is an
integral curve in ${\mathbb P}^3$ of degree~$4$, arithmetic genus~$1$, with one
node.  This is a complete intersection in ${\mathbb P}^3$, hence $S$ is a
complete intersection of two quadric hypersurfaces in
${\mathbb P}^4$ \cite[$1.3.3$]{M}, so it must contain every trisecant of
$C$.  But $S$ contains no lines except $L$, so $C$ has a unique trisecant. 
Since $C$ is linked to a double structure on $L$, $C$ is in the $CI$-liaison
class of two skew lines, so $C$ is in ${\mathcal L}_1$.  Note that $C-H$ is not
effective on $S$, so $C$ cannot be obtained by an elementary Gorenstein
biliaison on
$S$.  However, it seems likely that $C$ also lies on a normal singular Del
Pezzo surface on which it can be obtained by an elementary Gorenstein biliaison
from two skew lines.

Thus we see that any smooth nondegenerate $(6,1)$ curve in ${\mathbb P}^4$ is
in ${\mathcal L}_1$.  The family of all such curves in ${\mathbb P}^4$ is
irreducible \cite{Ein}.  The general type with two trisecants (Case 2a) is
obtained by an elementary Gorenstein biliaison from a curve of degree $2$ in
${\mathcal L}_0$, while the special type (Case~1) with infinitely many
trisecants is obtained by Gorenstein biliaison from a curve of degree $3$ in
${\mathcal L}_0$.

\bigskip
\noindent
{\bf Example 4.5.}  We consider nonsingular degree $7$ genus $2$ curves in
${\mathbb P}^4$, not contained in any ${\mathbb P}^3$.  The family $H_{7,2}$ of
all of these curves is irreducible, by Ein \cite{Ein}.  We see $h^1({\mathcal
I}_C(1)) = 1$, and there exist such curves with Rao module $k$ on a Del Pezzo
surface (see below), so by semicontinuity, the general such curve has Rao
module $k$, i.e., it is in ${\mathcal M}_1$.

Next, note that $h^0({\mathcal I}_C(2)) \ge 2$.  If $h^0({\mathcal I}_C(2)) >
2$, then the intersection of three quadric surfaces would either be a curve of
degree $8$, and then $C$ would be linked to a line, hence ACM, which is
impossible; or it would be a surface of degree $3$, but there are no $(7,2)$
curves on surfaces of degree $3$ in ${\mathbb P}^4$.  Hence $h^0({\mathcal
I}_C(2)) = 2$, and $h^1({\mathcal I}_C(2)) = 0$, so by Castelnuovo--Mumford
regularity, $C \in {\mathcal M}_1$.  Thus all curves of $H_{7,2}$ are in
${\mathcal M}_1$.

Now look on a Del Pezzo surface $S$, and let $C = (4;2,1^3,0)$.  Then $C$ is a
smooth $(7,2)$ curve, and $C - H = (1;1,0^3,-1)$ is a disjoint union of a conic
and a line.  Thus $C \in {\mathcal L}_1$, and $C$ is obtained by an elementary
Gorenstein biliaison from a curve of degree $3$ in ${\mathcal L}_0$.  Note that
$C$ has exactly four mutually skew trisecants, namely the lines $F_{25}$,
$F_{35}$, $F_{45}$, and $G$ on $S$.

If $C$ is any smooth $(7,2)$ curve, we have seen that $h^0({\mathcal I}_C(2)) =
2$.  Let $S$ be the complete intersection surface $F_2 \cdot F'_2$ of two
quadric surfaces containing $C$.  Then $S$ is uniquely determined by $C$.  It
is a surface of degree $4$, with sectional genus $1$, but it may be singular. 
However, it must be irreducible, and hence has at most a line of singular
points.  Therefore $C$ is an almost Cartier divisor on $S$, and we can apply
the theory of generalized divisors.  There is an exact sequence
\[
0 \rightarrow {\mathcal O}_S \rightarrow {\mathcal L}(C) \rightarrow
\omega_C(1) \rightarrow 0
\]
\cite[$2.1$]{GD}, making use of the fact that $\omega_S = {\mathcal O}_S(-1)$. 
Twisting by $-1$, and taking cohomology, we obtain
\[
0 \rightarrow H^0({\mathcal O}_S(-1)) \rightarrow H^0({\mathcal L}(C-H))
\rightarrow H^0(\omega_C) \rightarrow H^1({\mathcal O}_S(-1)),
\]
where $H$ denotes the hyperplane class on $S$.  The two outside groups are $0$,
and $H^0(\omega_C)$ has dimension $2$, so $H^0({\mathcal L}(C-H)) \ne 0$.  This
shows that $C-H$ is effective on $S$.  It is a divisor of degree $3$, and must
be in ${\mathcal L}_0$, so we see that any $C$ in ${\mathcal M}_1$ is in
${\mathcal L}_1$, and is obtained by an elementary Gorenstein biliaison from a
curve of degree $3$ in ${\mathcal L}_0$.

For an example of a special $(7,2)$ curve, let $S_0 = {\mathbb P}^1 \times
{\mathbb P}^1$.  Let $\Gamma$ be a line of bidegree $(1,0)$, and fix an
involution $\sigma$ on $\Gamma$.  Take $\vartheta$ to be the linear system of
those curves of bidegree $(1,2)$ on $S_0$ meeting $\Gamma$ in a pair of the
involution
$\sigma$.  Then $\vartheta$ maps $S_0$ to a surface $S$ of degree $4$ in
${\mathbb P}^4$ with a double line $L_0$ (the image of $\Gamma$).  If $C_0$ is
a general curve of bidegree $(2,3)$ on $S_0$, then the image $C$ of $C_0$ in
$S$ is a smooth $(7,2)$ curve meeting $L_0$ in three points.  It has four
trisecants, namely the double line $L_0$ and the three rulings (images of
$(0,1)$ curves in $S_0$) that meet $L_0$ at the points where $C$ meets $L_0$. 
This curve is different from the general ones described above, because three of
the trisecants meet the fourth one.  Because of the general result above, $C$
must arise by an elementary Gorenstein biliaison on $S$, but in this case the
curve of degree $3$ in ${\mathcal L}_0$ will be a nonreduced curve containing a
double structure on the line $L_0$.

Next we look at a general Castelnuovo surface $S'$.  On this surface, there are
three different kinds of smooth $(7,2)$ curves, distinguished by their
self-intersections, namely
\[
\begin{array}{rll}
C_1 &= &(4;2,1^5,0^2) \\
C_2 &= &(5;2^4,1^3,0) \\
C_3 &= &(5;1^4,2^4)
\end{array} \hskip 0.5 in 
\begin{array}{rll}
C_1^2 &= &7 \\
C_2^2 &= &6 \\
C_3^2 &= &5.
\end{array}
\]
Of these $C_1$ is obtained by an elementary Gorenstein biliaison on $S'$ from
two skew lines, while $C_2 - H$ and $C_3 - H$ are not effective.  Since we have
seen above that every smooth $(7,2)$ curve arises by elementary Gorenstein
biliaisons from a degree $3$ curve in ${\mathcal L}_0$, this gives examples of
curves that may be obtained by two different routes by elementary Gorenstein
biliaisons from curves of two different degrees in ${\mathcal L}_0$.

In fact, I claim that every general $(7,2)$ curve arises also as a curve of
type $C_1$ on a smooth Castelnuovo surface.  To prove this in detail is rather
long, so I will just give a sketch.  Start with a smooth $(7,2)$ curve $C$ on a
smooth Del Pezzo surface $S$, say $C = (4;2,1^3,0)$ as before.  Choose a
twisted cubic curve $D$ and a conic $\Gamma$ so that $C + D + \Gamma = 3H$. 
(For example $D = (3;1^4,2)$ and $\Gamma = (2;0,1^4)$.)  Let $\Pi$ be the plane
containing $\Gamma$.  Then the two quadric hypersurfaces containing $C$ meet
$\Pi$ in $\Gamma$, so a linear combination of them contains $\Pi$.  So we may
assume $S = F_2 \cdot F'_2$ where $F_2$ contains $\Pi$.  By construction,
there are cubic hypersurfaces $F_3$ containing $C+D+\Gamma$.  Such an $F_3$
will meet $\Pi$ in $\Gamma$ plus a line.  Adjusting $F_3$ by a linear form
times $F'_2$, we may assume that $F_3$ contains $\Pi$.  Now $F_2 \cdot F_3 =
\Pi \cup S'$, where $S'$ is an ACM surface of degree $5$, hence a Castelnuovo
surface.  Now one can verify that $C$ on $S'$ is a curve with self-intersection
$7$, like $C_1$ above, and that $C_1-H$ is effective and represented by a curve
of degree $2$ in ${\mathcal L}_0$.

In conclusion, we see that every smooth $(7,2)$ curve is in ${\mathcal L}_1$,
and can be obtained by elementary Gorenstein biliaison from ${\mathcal L}_0$,
in general by two different routes.  This is in contrast to the $(6,1)$ case
above, where the curves are divided into two types, distinguished by which
component of ${\mathcal L}_0$ they arise from.

\bigskip
\noindent
{\bf Example 4.6.}  We consider smooth $(10,6)$ curves in ${\mathbb P}^4$ (cf.\
Example $3.8$ above).

By Riemann--Roch applied to ${\mathcal O}_C(1)$ we see that an ACM $(10,6)$
curve is nonspecial.  Also we see that ${\mathcal O}_C(1)$ is special if and
only if ${\mathcal O}_C(1)$ is a canonical divisor, and this is equivalent to
$h^1({\mathcal I}_C(1)) = 1$.  The $(10,6)$ curves in ${\mathbb P}^4$ with
${\mathcal O}_C(1)$ special are all projections of the canonical curves of
genus $6$ in ${\mathbb P}^5$.  Since these form an irreducible family, we see
that their projections, the canonical $(10,6)$ curves in ${\mathbb P}^4$, form
an irreducible family, and they all have $h^1({\mathcal I}_C(1)) = 1$.  A
general such curve has Rao module $k$ in degree $1$, i.e., it is in ${\mathcal
M}_1$.  To see this, by semicontinuity, it is sufficient to exhibit one such. 
Let $C = (7;2^9,0)$ on a general Bordiga surface.  This curve is the image of a
plane septic curve with $9$ double points, embedded in ${\mathbb P}^4$ by the
linear system of quartics passing through the double points (and one further
point).  These are adjoint curves to $C$ and so cut out the canonical linear
series.  Hence $C$ is a canonical curve.  Now $C-H = (3;1^9,-1)$ is the
disjoint union of a plane cubic curve and a line, which is a degree $4$ curve
in ${\mathcal L}_0$.  Hence $C \in {\mathcal L}_1$ is obtained by an elementary
Gorenstein biliaison from ${\mathcal L}_0$.

Now let us consider $(10,6)$ curves in ${\mathcal M}_2$, i.e., with Rao module
$k$ in degree $2$.  Examples of such can be found on a general Castelnuovo
surface $S$, for example $C_1 = (6;3,2,1^6)$ and $C_2 = (6;2^4,1^4)$.  Note
that the curves $C_1$ are trigonal, while the curves $C_2$ can have general
moduli.  Both types can be obtained by Gorenstein biliaison on $S$, since $C_1
- H = (2;1^2,0^6)$ and $C_2 - H = (2;0,1^3,0^4)$ are both $(5,0)$ curves, hence
in ${\mathcal L}_1$.

An examination of curves of minimal genus in ${\mathcal L}_1$ and ACM curves of
minimal genus in ${\mathbb P}^3$ shows that the only way to obtain a $(10,6)$
curve in ${\mathcal M}_2$ by Gorenstein biliaison is from a $(5,0)$ curve in
${\mathcal L}_1$ on a surface of degree $5$ and sectional genus $2$ in
${\mathbb P}^4$, like the Castelnuovo surfaces.

Another example of a $(10,6)$ curve in ${\mathcal M}_2$ is obtained by the
curve $C$ formed by the intersection of a smooth quintic elliptic scroll $V$
with a hypersurface $F$ of degree $2$.  If we take $F$ to be a smooth quadric
hypersurface, then by Klein's theorem \cite[II.Ex.~$6.5$d]{AG} it contains no
surfaces of odd degree, so $C$ cannot be obtained by ascending Gorenstein
biliaison from a curve in ${\mathcal L}_1$.  

However one can show that $C$ is in the $G$-liaison class of two skew lines by
the following method, suggested by the referee.  First note that two general
cubic hypersurfaces $F_3,F'_3$ containing $V$ will link $V$ to be a Veronese
surface $W$ in ${\mathbb P}^4$.  Thus $C$ is linked by the complete
intersection $F_2 \cap F_3 \cap F'_3$ to a curve $C' \subseteq W$, which is $W
\cap F_2$.  The curve $C'$ is an $(8,3)$ curve, obtained from a plane curve of
degree~$4$ by the $2$-uple embedding of ${\mathbb P}^2$ and projection to
${\mathbb P}^4$.

Now $W$ is not an ACM surface, but if we take a hyperplane section $\Gamma = W
\cap {\mathbb P}^3$, then $\Gamma$ is a $(4,0)$ curve in ${\mathbb P}^3$.  It
is contained in a unique nonsingular quadric surface $Q \subseteq {\mathbb
P}^3$, and the union $W \cup Q$, meeting along $\Gamma$, will be an ACM surface
of degree~$6$ in ${\mathbb P}^4$.  We regard $C' \subseteq W$ as a curve on the
surface $W \cup Q$.  Now one can show (I leave some details to the reader) that
$2H-K-C'$ on the surface $W \cup Q$ (where $H,K$ denote the hyperplane section
and canonical divisor) is a curve $D \cup \Gamma'$, where $\Gamma'$ is a
$(4,0)$ curve in ${\mathbb P}^3$, and $D$ is a conic, not in ${\mathbb P}^3$,
meeting $\Gamma'$ in two points.  Since $2H-K$ is an arithmetically Gorenstein
curve on $W \cup Q$, we have thus linked $C$ to $C'$ and then to $D \cup
\Gamma'$.

For the last step, we take a quadric surface $Q'$ containing $D$ and meeting
the quadric $Q$ (which contains $\Gamma'$) in a conic.  Then $Q \cup Q'$ is a
complete intersection quartic surface in ${\mathbb P}^4$, and on $Q\cup Q'$, $D
\cup \Gamma' - H$ is $2$ skew lines.

Thus $C$ is an example of a curve with Rao module $k$, that cannot be obtained
by ascending Gorenstein biliaison from a minimal curve, and yet is in the
$G$-liaison class of $2$ skew lines.

\bigskip
\noindent
{\bf Example 4.7.} For our last example, we consider smooth $(11,7)$ curves in
${\mathcal M}_2$.  To construct such curves on a general Bordiga surface $S$,
take $C = (6;2^3,1^7)$. 
This is an $(11,7)$ curve, and
$C-H = (2;1^3,0^7)$ is a smooth $(5,0)$ curve on $S$.  Since $(5,0)
\in {\mathcal L}_1$ by $(4.3)$ above, we see that $C \in {\mathcal L}_2$, and
is obtained from a minimal curve by two elementary Gorenstein biliaisons.

Next, I claim the only way to obtain an $(11,7)$ curve in ${\mathcal L}_2$ by
two elementary $G$-biliaisons is the one just described.  Indeed, the curves of
minimal genus in ${\mathcal L}_1$ of degrees $4$ to $7$ are $(4,0)$, $(5,0)$,
$(6,1)$, $(7,2)$.  The minimal genus of ACM curves in ${\mathbb P}^3$ of
complementary degree are $(7,5)$, $(6,3)$, $(5,2)$, $(4,1)$, which will give
rise respectively to curves $(11,8)$, $(11,7)$, $(11,8)$, $(11,9)$ in
${\mathcal L}_2$.  So an $(11,7)$ curve obtained by elementary $G$-biliaisons
must be on the Bordiga surface or its specialization.

Since $h^0({\mathcal O}_C(1)) = 5$, we see that ${\mathcal O}_C(1)$ is nonspecial, so we
can compute the dimension of the Hilbert scheme of $(11,7)$ curves in ${\mathbb
P}^4$ (which is irreducible by Ein \cite{Ein}), by the usual formula $5d+1-g$. 
Thus the Hilbert scheme has dimension $49$.

Now let us count the curves obtained by the construction above.  The Bordiga
surface moves in a family of dimension $36$ (use, for example, the formula of
Ellingsrud \cite{E}).  To find the dimension of the linear system $|C|$ on $S$,
we use the exact sequence
\[
0 \rightarrow {\mathcal O}_S \rightarrow {\mathcal O}_S(C) \rightarrow {\mathcal O}_C(C)
\rightarrow 0.
\]
Thus $\dim_S |C| = h^0({\mathcal O}_C(C))$.  For the curve $C$ of type
$(6;2^3,1^7)$ mentioned above, we find $C^2 = 17$, so the divisor $C^2$ is
nonspecial on $C$, and by Riemann--Roch, $h^0({\mathcal O}_C(C)) = C^2 + 1 - g
= 11$.  Thus the dimension of the family of all curves of this type on Bordiga
surfaces is $\le 11 + 36 = 47$.  In particular, these curves are not general
among all $(11,7)$ curves.

But we wish to show more, namely that a general $(11,7)$ curve does not lie on
a Bordiga surface.  So suppose now that $C$ is any $(11,7)$ curve on a Bordiga
surface.  I claim $C^2 \le 17$.  Indeed, we have $2g - 2 = C^2 + C.K$.  On the
Bordiga surface let $C = (a;b_1,\dots,b_{10})$.  The canonical divisor $K$ can
be written $K = (1;0^{10}) - H$.  So $2g - 2 = C^2 + a - d$, and $C^2 = 2g - 2
+ d - a = 23 - a$.  But in order to get a curve of genus~$7$, we must have 
$a \ge 6$.  Thus $C^2 \le 17$.  Then the same argument as above shows that
$h^0({\mathcal O}_C(C)) \le 11$, and we get the same dimension count, unless
${\mathcal O}_C(C)$ is a special divisor.  But in that case $C^2 \le 12$, and
by Clifford's theorem $h^0({\mathcal O}_C(C)) \le 7$.  Thus a general $(11,7)$
curve does not lie on a Bordiga surface.

Next, observe that for any $(11,7)$ curve in ${\mathbb P}^4$, $h^0({\mathcal O}_C(2))
= 16$, so necessarily $h^1({\mathcal I}_C(2)) \ge 1$.  Since we have
constructed curves $C$ with $h^1({\mathcal I}_C(2)) = 1$ and the other
$h^1({\mathcal I}_C(n)) = 0$ for $n \ne 2$, we conclude by semicontinuity that
the general $(11,7)$ curve in ${\mathbb P}^4$ has Rao module $k$ in degree $2$,
i.e., it lies in ${\mathcal M}_2$.  Since the general such curve does not lie
on a Bordiga surface, by the above remarks, it cannot be obtained by ascending
Gorenstein biliaisons from ${\mathcal L}_0$.

It is conceivable that the general $(11,7)$ curve is linked by some ascending
and descending $G$-liaisons to two skew lines, but this seems unlikely, so we
propose the general $(11,7)$ curve in ${\mathbb P}^4$ as a possible curve with
Rao module $k$, not in the $G$-liaison class of two skew lines.

\section{Conclusion}

The examples presented in this paper would lead me to expect that for ACM
schemes of codimension $\ge 3$, some may be obtained by elementary Gorenstein
biliaisons from a scheme of degree one; a broader class may be obtained by
ascending and descending $G$-liaisons from a scheme of degree one; but that a
general ACM scheme of high degree may not be in the $G$-liaison class of a
complete intersection.  Example $3.9$ shows that at least one of the Questions
$1.3$a, $1.3$b has no for an answer.  Namely, for the general ACM $(20,26)$
curve in ${\mathbb P}^4$, we must have either

\begin{itemize}
\item[a)] it is ACM and not glicci, or
\item[b)] it is glicci, but cannot be obtained by ascending Gorenstein
biliaisons from a curve.
\end{itemize}

For curves in ${\mathbb P}^n$, $n\ge 4$, with a given Rao module $M$, I would
expect that the minimal curves form an infinite union of irreducible families;
some curves in the family may be obtained by a sequence of ascending elementary
Gorenstein biliaisons from a minimal curve; a larger class may be obtained by
ascending and descending Gorenstein liaisons from a minimal curve; but that a
general curve of high degree and genus with Rao module $M$ is not in the
$G$-liaison class of a minimal curve.  Example $4.6$ gives an example of a
smooth $(10,6)$ curve with Rao module $k$, that is in the Gorenstein liaison
class of a minimal curve, but cannot be obtained by ascending Gorenstein
biliaison from a minimal curve.  Example
$4.7$ shows that either Question
$1.4$ has no for an answer, or the deformation is necessary in Question
$1.5$b.  Indeed, for the general $(11,7)$ curve in ${\mathbb P}^4$ we must have
either

\begin{itemize}
\item[a)] it has Rao module $k$, but is not in the $G$-liaison class of two skew
lines, or
\item[b)] it is in the liaison class of two skew lines, but cannot be obtained
from a minimal curve by ascending Gorenstein biliaisons.
\end{itemize}

Based on this evidence, I would expect a no answer to Questions $1.3$, $1.4$,
$1.5b$.  I have no idea about Question $1.2$ since in this paper I used only
strict Gorenstein liaisons and biliaisons.  Also the question of even or odd
liaison has not been addressed here, since it is irrelevant for ACM schemes and
curves with Rao module $k$.  This is a question that merits further study.

\end{document}